%% file: higher.tex
\documentclass[11pt,a4paper]{amsart}

\usepackage{latexsym}
\usepackage{amsmath}
\usepackage{amsthm}
\usepackage{amssymb}
\usepackage{vmargin}
\usepackage{amscd}
\usepackage{stmaryrd}
\usepackage{euscript}
\usepackage{mathrsfs}
\usepackage{amscd}
\usepackage[all]{xy}
\usepackage{xr}

\DeclareMathAlphabet{\mathpzc}{OT1}{pzc}{m}{it}

\include{newthms}

\include{newcommands}

\begin{document}
\begin{abstract}
We generalise the techniques of \cite{ddt2} 
to describe  derived deformations in simplicial categories. This allows us to consider deformation problems with higher automorphisms, such as chain complexes (which have homotopies) and stacks (which have $2$-automorphisms). We also give a general approach for studying deformations of diagrams.
\end{abstract}

\title{Derived deformations of Artin stacks}
\author{J.P.Pridham}
\thanks{The author was supported during this research by Trinity College, Cambridge.}
\maketitle


\section*{Introduction}

This paper is  motivated by the wish  to describe derived deformations of an algebraic stack. In \cite{olssonstack} and \cite{aoki}, it was shown that deforming  an algebraic  stack can be regarded as a special case of deforming a simplicial algebraic space. The category of simplicial spaces has a natural simplicial structure (meaning that the $\Hom$-sets can be enriched to give simplicial sets), and the $2$-groupoid of deformations of an algebraic stack can be recovered from this simplicial structure.

After reviewing background material from \cite{ddt1} in \S \ref{one}, 
 we  introduce derived deformation complexes (DDCs) In Section \ref{ddcsn}; these extend the SDCs of \cite{paper2}  to simplicial categories. We then adapt the various constructions of \cite{ddt1}, showing how to associate derived deformation functors to DDCs, and how to compare them with derived deformation functors coming from SDCs.

Section \ref{constrsn} adapts the ideas of \cite{paper2}, showing how to associate DDCs to bialgebraic deformation problems in simplicial categories. In \S \ref{weakinvart}, we  show how  deformations of morphisms and diagrams can be used  to compare deformations of weakly equivalent objects.

Several simple examples of such problems are considered in Section \ref{egsn}: chain complexes (with more interesting variants in Remarks \ref{moreinterest}), simplicial complexes and simplicial algebras.

The motivating example of algebraic stacks is finally considered in Section \ref{stacksn}. We first describe derived deformations of simplicial affine schemes (\S \ref{calgsn}), then show in \S \ref{nicestacks} how to adapt this to describe derived deformations of an algebraic stack $\fX$, with an indication in Remark \ref{inftyrk} of  how this approach also works for Artin $n$-stacks. The idea is to consider derived deformations of  a suitable hypercovering $X_{\bt}$ of $\fX$. To see that this does, indeed, extend the $2$-groupoid of  deformations of $\fX$, we establish comparisons with Olsson's $\Ext$-groups of the cotangent complex (\S \ref{cfolsson}) and Aoki's description of the deformation $2$-groupoid (\S \ref{cfaoki}).

\tableofcontents

\section{Derived deformation functors}\label{one}

With the exception of \S \ref{quotsp},   the definitions and results in this section can all be found in \cite{ddt1}. 
Fix a complete local Noetherian ring $\L$, with maximal ideal $\mu$ and residue field $k$. 

\subsection{Simplicial Artinian rings}

\begin{definition}
Let $\C_{\L}$ denote the category of local Artinian $\L$-algebras with residue field $k$.
 We define $s\C_{\L}$ to be the category of Artinian simplicial  local $\L$-algebras, with residue field $k$. 
\end{definition}

\begin{definition}\label{N^s}
Given a simplicial complex $V_{\bt}$, recall that the normalised chain complex $N^s(V)_{\bt}$ is given by  $N^s(V)_n:=\bigcap_{i>0}\ker (\pd_i: V_n \to V_{n-1})$, with differential $\pd_0$. The simplicial Dold-Kan correspondence says that $N^s$ gives an equivalence of categories between simplicial complexes and non-negatively graded chain complexes in any abelian category. Where no ambiguity results, we will denote $N^s$ by $N$.
\end{definition}

\begin{lemma}\label{cotdef}
A simplicial complex $A_{\bt}$  of local $\L$-algebras with residue field $k$ and maximal ideal $\m(A)_{\bt}$   is Artinian if and only if:
\begin{enumerate}
\item the normalisation $N(\cot A)$ of the cotangent space $\cot A:=\m(A)/(\m(A)^2+\mu \m(A))$  is finite-dimensional (i.e. concentrated in finitely many degrees, and finite-dimensional in each degree). 
\item For some $n>0$, $\m(A)^n=0$.
\end{enumerate} 
\end{lemma}
\begin{proof}
\cite{ddt1} Lemma \ref{ddt1-cotdef}
\end{proof}

 As in \cite{descent},  we say that a functor is left exact if it preserves all finite limits. This is equivalent to saying that it preserves final objects and fibre products.

\begin{definition}\label{spdef}
Define $\Sp$ to be  the category  of left-exact functors from $\C_{\L}$ to $\Set$. 
Define  $c\Sp$ to be the category of left-exact functors from $s\C_{\L}$ to $\Set$.
\end{definition}

\begin{definition}
Given a functor $F:\C_{\L} \to \Set$, we write $F:s\C_{\L} \to \Set$ to mean $A \mapsto F(A_0)$  (corresponding to the inclusion $\Sp \into c\Sp$).
\end{definition}

\subsection{Properties of morphisms}

\begin{definition}\label{smoothdef}
As in \cite{Man}, we say that a functor $F:\C_{\L}\to \Set$ is smooth if for all surjections $A \to B$ in $\C_{\L}$, the map $F(A) \to F(B)$ is surjective. 
\end{definition}

\begin{definition}
We say that a map $f:A \to B$ in $s\hat{\C}_{\L}$ is acyclic if $\pi_i(f):\pi_i(A) \to \pi_i(B)$ is an isomorphism of pro-Artinian $\L$-modules for all $i$.  $f$ is said to be surjective if each $f_n:A_n \to B_n$ is  surjective.
\end{definition}

Note that for any simplicial abelian group $A$, the homotopy groups can be calculated by $\pi_iA \cong \H_i(NA)$, the homology groups of the normalised chain complex. These in turn are isomorphic to the homology groups of the unnormalised chain complex associated to $A$. 

\begin{definition}
We define a small extension $e:I \to A \to B$ in $s\C_{\L}$ to consist of a surjection $A \to B$ in $s\C_{\L}$ with kernel $I$, such that $\m(A)\cdot I=0$. Note that this implies that $I$ is a simplicial complex of $k$-vector spaces.
\end{definition}

\begin{lemma}\label{small}
Every surjection in $s\C_{\L}$ can be factorised as a composition of small extensions. Every acyclic surjection in $s\C_{\L}$ can be factorised as a composition of acyclic small extensions.  
\end{lemma}
\begin{proof}
\cite{ddt1} Lemma \ref{ddt1-small}.
\end{proof}

\begin{definition}
We say that a morphism $\alpha:F\to G$ in $c\Sp$  is smooth if for all small extensions $A \onto B$ in $s\C_{\L}$, the map $F(A) \to F(B)\by_{G(B)}G(A)$ is surjective. 

Similarly, we call $\alpha$ quasi-smooth if for all acyclic small extensions $A \to B$ in $s\C_{\L}$, the map $F(A) \to F(B)\by_{G(B)}G(A)$ is surjective.
\end{definition}

\begin{lemma}\label{ctod} A morphism  $\alpha:F\to G$ in $\Sp$ is smooth if and only if the induced morphism between the  objects $F, G \in c\Sp$ is quasi-smooth, if and only if it is smooth.
\end{lemma}
\begin{proof}
\cite{ddt1} Lemma \ref{ddt1-ctod}.
\end{proof}

\subsection{Derived deformation functors}

\begin{definition}
Define the  $sc\Sp$ to be the category of left-exact functors from $s\C_{\L}$ to the category $\bS$ of simplicial sets.
This is equivalent to the category of simplicial cosimplicial objects in $\Sp$.

Define $s\Sp$ to be  the category  of left-exact functors from $\C_{\L}$ to $\bS$.
\end{definition}

\begin{definition}\label{scspqsdef}
 A morphism $\alpha:F\to G$ in    $sc\Sp $ is  said to be smooth if 
\begin{enumerate}
\item[(S1)]
for every acyclic surjection $A \to B$ in $s\C_{\L}$, the map $F(A)\to F(B)\by_{G(B)}G(A)$ is a trivial fibration in $\bS$; 
\item[(S2)]
for every surjection $A \to B$ in $s\C_{\L}$, the map $F(A)\to F(B)\by_{G(B)}G(A)$ is a surjective fibration in $\bS$.
\end{enumerate}

A morphism $\alpha:F\to G$ in    $sc\Sp $  is  said to be quasi-smooth if it satisfies (S1) and
\begin{enumerate}
\item[(Q2)]
for every surjection $A \to B$ in $s\C_{\L}$, the map $F(A)\to F(B)\by_{G(B)}G(A)$ is a  fibration in $\bS$.
\end{enumerate}
\end{definition}

\begin{definition}
Given $A \in s\C_{\L}$ and a finite simplicial set $K$, define $A^K \in \C_{\L}$ by 
$$
(A^K)_i:=\Hom_{\bS}(K\by \Delta^i, A)\by_{\Hom_{\Set}(\pi_0K, k)}k.
$$
\end{definition}

\begin{definition}\label{underline}
Given   $F \in sc\Sp$, define $\underline{F}:s\C_{\L}\to \bS$ by 
$$
\underline{F}(A)_n:= F_n(A^{\Delta^n}).
$$
 
For $F \in c\Sp$, we may regard $F$ as an object of $sc\Sp$ (with the constant simplicial structure), and then define $\underline{F}$  as above.
\end{definition}

\begin{lemma}\label{settotop} A map $\alpha:F\to G $ in $c\Sp$ is smooth (resp. quasi-smooth) if and only if the induced map of functors  $\underline{\alpha}:\underline{F}\to \underline{G}$  is smooth (resp. quasi-smooth) in $sc\Sp$.
\end{lemma}
\begin{proof}
\cite{ddt1} Lemma \ref{ddt1-settotop}.
\end{proof}

The following Lemma will provide many examples of functors which are quasi-smooth but not smooth. 
\begin{lemma}\label{sm7}
If $F\to G$ is a quasi-smooth map of functors $F,G:s\C_{\L} \to \bS$, and $K \to L$ is a cofibration in $\bS$, then 
$$
F^L \to F^K\by_{G^K}G^L
$$
is quasi-smooth.
\end{lemma}
\begin{proof}
This is an immediate consequence of the fact that $\bS$ is a simplicial model category, following from axiom SM7, as given in \cite{sht} \S II.3.
\end{proof}

The following lemma is a consequence of  standard properties of fibrations and trivial fibrations in $\bS$.
\begin{lemma}\label{basechange}
If $F\to G$ is a quasi-smooth map of functors $F,G:s\C_{\L} \to \bS$, and $H \to G$ is any map of functors, then $F\by_GH \to H$ is quasi-smooth. 
\end{lemma}

\begin{definition}\label{pioqs}
A map $ \alpha:F\to G$ of functors   $F,G:\C_{\L}\to \bS$ is said to be smooth (resp. quasi-smooth, resp. trivially smooth)  if for all surjections $A \onto B$ in $\C_{\L}$, the maps 
$$
 F(A) \to F(B)\by_{G(B)}G(A)
$$ 
  are surjective fibrations (resp. fibrations, resp. trivial fibrations).
\end{definition}

\begin{proposition}\label{smoothchar}
A map  $ \alpha:F\to G$  of left-exact functors   $F,G:\C_{\L}\to \bS$ is smooth if and only if the maps  $F_n\xra{\alpha_n} G_n$ of functors   $F_n,G_n:\C_{\L}\to \Set$ are all smooth.
\end{proposition}
\begin{proof}
\cite{ddt1} Proposition \ref{ddt1-smoothchar}.
\end{proof}

\begin{proposition}\label{ctodtnew}
If a morphism $F\xra{\alpha} G$ of left-exact functors $F,G:s\C_{\L} \to \bS$ is such that the maps
$$
\theta: F(A)\to F(B)\by_{G(B)}G(A)
$$
 are surjective fibrations  for all acyclic small extensions $A \to B$, then $\underline{\alpha}:\underline{F} \to \underline{G}$ is quasi-smooth (resp. smooth) if and only if  $\theta$ is a fibration (resp. surjective fibration) for all small extensions $A \to B$.
\end{proposition}
\begin{proof}
\cite{ddt1} Proposition \ref{ddt1-ctodtnew}.
\end{proof}

\begin{definition}\label{weakdef}
We will say that a morphism $\alpha: F \to G$ of quasi-smooth objects of $sc\Sp$ is a weak equivalence if, for all $A \in s\C_{\L}$, the maps $\pi_i F(A) \to \pi_iG(A)$ are isomorphisms for all $i$.
\end{definition}

\subsection{Quotient spaces}\label{quotsp}

\begin{definition}\label{quotdefn}
Given  functors $X:s\C_{\L} \to \bS$ and $G: s\C_{\L} \to s\Gp$, together with a right action of $G$ on $X$, define the quotient space by
$$
[X/G]_n= (X\by^{G}WG)_n=  X_n\by G_{n-1}\by G_{n-2}\by \ldots G_0,
$$
with operations as standard for universal bundles (see \cite{sht} Ch. V). Explicitly:
\begin{eqnarray*}
\pd_i(x,g_{n-1},g_{n-2}, \ldots,g_0)&=& \left\{ \begin{matrix}(\pd_0x*g_{n-1}, g_{n-2}, \ldots, g_0)& i=0;\\
 (\pd_ix, \pd_{i-1}g_{n-1}, \ldots, (\pd_0g_{n-i})g_{n-i-1}, g_{n-i-2}, \ldots, g_0) & 0<i<n;\\
(\pd_nx, \pd_{n-1}g_{n-1},\ldots, \pd_1g_1) & i=n;
\end{matrix}\right.\\
\sigma_i(x,g_{n-1},g_{n-2}, \ldots,g_0)&=& (\sigma_ix, \sigma_{i-1}g_{n-1}, \ldots, \sigma_0g_{n-i},e,g_{n-i-1}, g_{n-i-2}, \ldots, g_0).
\end{eqnarray*}
The space $[\bullet/G]$ is also denoted $\bar{W}G$, and is a model for the classifying space $BG$ of $G$. 
\end{definition}

\begin{lemma}
If $G:s\C_{\L} \to s\Gp$ is smooth, then  $\bar{W}G$ is smooth.
\end{lemma}
\begin{proof}
For any surjection $A \to B$, we have $G(A) \to G(B)$ fibrant and surjective on $\pi_0$, which by \cite{sht} Corollary V.6.9 implies that $\bar{W}G(A) \to \bar{W}G(B)$ is a fibration. If $A \to B$ is also acyclic, then everything is trivial by properties of $\bar{W}$ and $G$.
\end{proof}

\begin{remark}
Observe that this is our first example of a quasi-smooth functor  which is not a right  Quillen functor for the simplicial model structure. The definitions of smoothness and quasi-smoothness  were designed with $\bar{W}G$ in mind.
\end{remark}

\begin{lemma}
If $X$ is quasi-smooth, then so is $[X/G] \to \bar{W}G$.
\end{lemma}
\begin{proof}
This follows from the observation that for any fibration (resp. trivial fibration) $Z \to Y$ of $G$-spaces, $[Z/G]\to [Y/G]$ is a fibration (resp. trivial fibration).
\end{proof}

\begin{corollary}\label{fibquot}
If $X$ is quasi-smooth and $G$  smooth, then  $[X/G]$ is quasi-smooth.
\end{corollary}
\begin{proof}
Consider the fibration $X \to [X/G] \to \bar{W}G$.
\end{proof}

\subsection{Cohomology and obstructions}\label{cohomology}

Given a quasi-smooth morphism  $\alpha:F\to G$  in $sc\Sp$, there exist $k$-vector spaces  $\H^i(F/G)$ for all $i \in \Z$. 

By \cite{ddt1} Corollary \ref{ddt1-cohowelldfn}, these have the
 property that for any  simplicial  $k$-vector space $V$ with finite-dimensional normalisation,
$$
\pi_m(F(k\oplus V)\by_{G(k\oplus V)}\{0\}) \cong \H^{-m}(F/G\ten V),
$$
where $V^2=0$ and 
$$
H^i(F/G\ten V):=\bigoplus_{n \ge 0} \H^{i+n}(F/G) \ten \pi_n(V).
$$

If $G=\bt$ (the one-point set), we write $\H^j(F):= \H^j(F/\bt)$.

We now have the following characterisation of obstruction theory:

\begin{theorem}\label{robs}
If  $\alpha:F \to G$ in $sc\Sp$ is quasi-smooth, then for any small extension $e:I \to A \xra{f} B$ in $s\C_{\L}$, there is a sequence of sets
$$
\pi_0(FA)\xra{f_*} \pi_0(FB\by_{GB}GA) \xra{o_e}  \H^1(F/G\ten I)
$$  
exact in the sense that the fibre of $o_e$ over $0$ is the image of $f_*$. Moreover,  there is a group action of $\H^0(F/G \ten I)$ on $\pi_0(FA)$ whose orbits are precisely the fibres of $f_*$. 

For any $y \in F_0A$, with $x=f_*y$, the fibre of $FA \to FB\by_{GB}GA$ over $x$ is isomorphic to $\ker(\alpha: FI \to GI)$, and the sequence above 
extends to a long exact sequence
$$\xymatrix@R=0ex{
\cdots  \ar[r]^-{f_*}&\pi_n(FB\by_{GB}GA,x) \ar[r]^-{o_e}& \H^{1-n}(F/G \ten I) \ar[r]^-{\pd_e} &\pi_{n-1}(FA,y)\ar[r]^-{f_*}&\cdots\\ 
\cdots \ar[r]^-{f_*}&\pi_1(FB\by_{GB}GA,x) \ar[r]^-{o_e}& \H^0(F/G \ten I)  \ar[r]^-{-*y} &\pi_0(FA).
}
$$
\end{theorem}
\begin{proof}
\cite{ddt1} Theorem \ref{ddt1-robs}.
\end{proof}

\begin{corollary}\label{weak} 
A map $\alpha:F \to G$ of quasi-smooth  $F,G\in sc\Sp$ is a weak equivalence if and only if 
the maps $\H^j(\alpha):\H^j(F) \to\H^j(G)$ are all isomorphisms.
\end{corollary}

\begin{corollary}\label{cohosmoothchar}
If $\alpha:F \to G$ is  quasi-smooth in $sc\Sp$, then $\alpha$ is smooth if and only if $\H^i(F/G)=0$ for all $i>0$.
\end{corollary}

\begin{proposition}\label{longexact}
Let $X,Y, Z:s\C_{\L}\to \bS$ be left-exact functors, with  $X \xra{\alpha} Y$ and $Y \xra{\beta} Z$ quasi-smooth.  There is then a long exact sequence
$$
\ldots \xra{\pd} \H^j(X/Y) \to \H^j(X/Z) \to \H^j(Y/Z) \xra{\pd} \H^{j+1}(X/Y) \to \H^{j+1}(X/Z) \to \ldots
$$
\end{proposition}
\begin{proof}
\cite{ddt1} Proposition \ref{ddt1-longexact}.
\end{proof}

\subsection{Model structures}\label{model}

\begin{theorem}\label{scspmodel}
There is a simplicial model structure on $sc\Sp$, for which the   fibrations are quasi-smooth morphisms, and weak equivalences between quasi-smooth objects are those given in Definition \ref{weakdef}.
\end{theorem}
\begin{proof}
This is \cite{ddt1} Theorem \ref{ddt1-scspmodel}.
\end{proof}

Thus the homotopy category $\Ho(sc\Sp)$ is equivalent to the category of quasi-smooth objects in $sc\Sp$, localised at the weak equivalences of Definition \ref{weakdef}.  

\begin{definition}\label{generalcoho}
Given any morphism $f:X \to Z$, we define $\H^n(X/Z):= \H^n(\hat{X}/Z)$, for $X \xra{i} \hat{X} \xra{p} Z$ a factorisation of $f$ with $i$ a geometric trivial cofibration, and $p$ a geometric fibration. 
\end{definition}

\subsubsection{Homotopy representability}

\begin{definition}\label{schless}
Define the category $\cS$ to consist of functors $F: s\C_{\L}\to \bS$ satisfying the following conditions:
\begin{enumerate}

\item[(A0)] $F(k)$ is contractible.

\item[(A1)] For all small extensions $A \onto B$ in $s\C_{\L}$, and maps $C \to B$ in $s\C_{\L}$, the 
map 
$F(A\by_BC) \to F(A)\by_{F(B)}^hF(C)$ is a  
weak equivalence, where $\by^h$ denotes homotopy fibre product.

\item[(A2)] For all acyclic small extensions $A \onto B$ in $s\C_{\L}$, the map $F(A) \to F(B)$ is a weak equivalence.
\end{enumerate}

Say that a natural transformation $\eta:F \to G$ between such functors is a weak equivalence if the maps $F(A) \to G(A)$ are weak equivalences for all $A\in s\C_{\L}$, and let $\Ho(\cS)$ be the category obtained by formally inverting all weak equivalences in $\cS$.
\end{definition}

\begin{theorem}\label{schrep}
There is a canonical equivalence between the geometric homotopy category $\Ho(sc\Sp)$ and the category $\Ho(\cS)$.
\end{theorem}
\begin{proof}
This is \cite{ddt1} Theorem \ref{ddt1-schrep}.
\end{proof}

\subsubsection{Equivalent formulations}

If $k$ is a field of characteristic $0$, then we may work with dg algebras rather than simplicial algebras.

\begin{definition}
Define $dg\C_{\L}$ to be the category of Artinian local differential $\N_0$-graded  graded-commutative  $\L$-algebras with residue field $k$.
\end{definition}

\begin{definition}
Define a map $A \to B$ in $dg\C_{\L}$ to be a small extension if it is surjective and the kernel $I$ satisfies $I\cdot \m(A)=0$. 
\end{definition}

\begin{definition}\label{sdgspdef}
Define $sDG\Sp$ to be the category of left exact functors from $dg\C_{\L}$ to $\bS$.
\end{definition}

\begin{definition}
Say a map $X \to Y$ in $sDG\Sp$ is quasi-smooth if for all small extensions $f:A \to B$ in $dg\C_{\L}$, the morphism
$$
X(A)\to Y(A)\by_{Y(B)}X(B)
$$ 
is a fibration in $\bS$, which is moreover a  trivial fibration if $f$ is acyclic. 
\end{definition}

\begin{definition}\label{weakdef2}
We will say that a morphism $\alpha: F \to G$ of quasi-smooth objects of $sDG\Sp$ is a weak equivalence if, for all $A \in s\C_{\L}$, the maps $\pi_i F(A) \to \pi_iG(A)$ are isomorphisms for all $i$.
\end{definition}

\begin{proposition}\label{sdgmodel}
There is a model structure on $sDG\Sp$, for which the   fibrations are quasi-smooth morphisms, and weak equivalences between quasi-smooth objects are those given in Definition \ref{weakdef2}.
 \end{proposition}
\begin{proof}
This is \cite{ddt1} Proposition \ref{ddt1-sdgmodel}.
\end{proof}

Most of the constructions from $s\C_{\L}$ carry over to $dg\C_{\L}$. However, there is no straightforward analogue of Definition \ref{underline}.

\begin{definition}\label{qrat}
Define the normalisation functor $N: s\C_{\L} \to dg\C_{\L}$ by mapping $A$ to its associated normalised complex $NA$, equipped with the Eilenberg-Zilber shuffle product (as in \cite{QRat}).
\end{definition}

\begin{definition}
Define $\Spf N^*:  sDG\Sp \to sc\Sp$ by mapping $X: dg\C_{\L} \to \bS$ to the composition $X\circ N:s\C_{\L} \to \bS$. Note that this is well-defined, since $N$ is left exact.
\end{definition}

\begin{theorem}\label{nequiv}
$\Spf N^*:  sDG\Sp \to sc\Sp$ is a right Quillen equivalence. 
\end{theorem}
\begin{proof}
This is \cite{ddt1} Theorem \ref{ddt1-nequiv}.
\end{proof}

In particular, this means that $\Spf N^*$ maps quasi-smooth morphisms to quasi-smooth morphisms, and induces an equivalence $ \oR\Spf N^*: \Ho( sDG\Sp) \to \Ho(sc\Sp)$.

\section{Derived deformation complexes}\label{ddcsn}

\subsection{Definitions}

\begin{definition}
Define a pre-SDC to consist of homogeneous functors $E^n: \C_{\L} \to \Set$, for $n \in \N_0$, together with maps 
$$
\begin{matrix}
\pd^i:E^n \to E^{n+1} & 1\le i \le n\\
\sigma^i:E^{n}\to E^{n-1} &0 \le i <n,
\end{matrix}
$$
an associative product $*:E^m \by E^n \to E^{m+n}$, with identity $1: \bullet \to E^0$,  such that:
\begin{enumerate}
\item $\pd^j\pd^i=\pd^i\pd^{j-1}\quad i<j$.
\item $\sigma^j\sigma^i=\sigma^i\sigma^{j+1} \quad i \le j$.
\item 
$
\sigma^j\pd^i=\left\{\begin{matrix}
			\pd^i\sigma^{j-1} & i<j \\
			\id		& i=j,\,i=j+1 \\
			\pd^{i-1}\sigma^j & i >j+1
			\end{matrix} \right. .
$
\item $\pd^i(e)*f=\pd^i(e*f)$.
\item $e*\pd^i(f)=\pd^{i+m}(e*f)$, for $e \in E^m$.
\item $\sigma^i(e)*f=\sigma^i(e*f)$.
\item $e*\sigma^i(f)=\sigma^{i+m}(e*f)$, for $e \in E^m$.
\end{enumerate}
\end{definition}

\begin{remark}
Note that a  pre-SDC is an SDC (in the sense of \cite{paper2} if and only if the spaces $E^n$ are smooth for all $n$.
\end{remark}

\begin{definition}\label{ddcdef} 
Define a pre-derived deformation complex (pre-DDC) $E$ to be
a simplicial complex $E_{\bt}$ of pre-SDCs. 

Given $K \in \bS$, observe that $E^n_K:=\Hom_{\bS}(K,E^n)$ is  a pre-SDC. 
\end{definition}

\begin{remark}
If each $E_m$ is an SDC, then Lemma \ref{smoothchar} implies that for all $n$, $E^n:\C_{\L}\to \bS$ is smooth. For $K \in \bS$ contractible, this implies that $E_K$ is an SDC.
\end{remark}

\begin{definition}\label{tandef}
Given a left-exact functor $F: \C_{\L} \to \Set$, define the tangent space $t_F$ (or $t(F)$) by $t_F:= F(k[\eps]/(\eps^2))$. Since  $k[\eps]/(\eps^2)$ is an abelian group object in $\C_{\L}$, $t_F$ is an abelian group. The endomorphisms $\eps \mapsto \lambda \eps$ of $k[\eps]/(\eps^2)$ make $t_F$ into a  vector space over $k$. 

Given a morphism $\alpha: F\to G$ of such functors, define the relative tangent space $\tan(F/G):= \ker(\tan F \to \tan G)$.
\end{definition}

\begin{definition}
Given a morphism $f:E \to F$ of  pre-SDCs for which each $f^n: E^n \to F^n$ is smooth, we may define cohomology groups $\H^*(E/F)$ as cohomology of the cosimplicial complex $\CC^{\bt}(E/F)$ given by $\CC^n(E/F):= \tan(E/F)$,  with cosimplicial structure defined as in \cite{paper2} \S\ref{paper2-sdc}.  
\end{definition}

\begin{definition}
Given a morphism   $E\to F$ of pre-DDCs, levelwise smooth in the sense that  each $f^i_n: E^i_n \to F^i_n$ is smooth, observe that the cohomology groups $\H^i(E/F)$ are simplicial vector spaces, and denote the corresponding  normalised chain complexes by $N\H^i(E/F)$. 
\end{definition}

\begin{definition}\label{q12}
A morphism $f:E \to F$ of pre-DDCs is said to be quasi-smooth  if:
\begin{enumerate}

\item[Q1.] for all $n,i\ge 0$, $E_n^i \to E_{\pd\Delta^n}^i\by_{F_{\pd\Delta^n}^i}F_n^i$ is smooth, and

\item[Q2.] for all $i>0$, $\H^i(E/F)$ is a constant simplicial complex, or equivalently

\item[Q2'.] for all $n> 0, i>0$, $N_n\H^i(E/F)=0$.

\end{enumerate}

Say that a pre-DDC $E$ is a DDC if it is quasi-smooth, i.e.  if $E \to \bt$ is quasi-smooth.
\end{definition}

\begin{definition}
Given a  levelwise smooth morphism $f:E \to F$ of pre-DDCs,  define the tangent chain cochain complex by $N\CC^{\bt}_{\bt}(E/F)$. 
\end{definition}

\begin{definition}\label{q12vect}
Say  that a simplicial cosimplicial complex $V \in sc\Vect_k$ is  quasi-smooth 
if  $\H_n(  NV^i)=0$ for all $n,i\ge 0$ and $\H^i(NV)_n=0$ for all $i>0$ and $n>0$.
\end{definition}

\begin{definition}\label{corner}
Given $V \in sc\Vect_k$ quasi-smooth, define a cochain complex $\lrcorner V$  by:
$$
(\lrcorner V)^n:= \left\{ \begin{matrix} V_0^n & n\ge 0\\ \H^0(N_{-n}V) & n <0, \end{matrix} \right. 
$$
with differential $d_c$ in non-negative degrees, and $d^s$ in negative degrees.

Given a  levelwise smooth morphism $f:E \to F$ of pre-DDCs, 
define the cohomology groups $\H^*(\lrcorner E/F):= \H^*(\lrcorner \CC^{\bt}(E/F))$, noting that these are given by
$$
\H^i(\lrcorner E/F)\cong \left\{\begin{matrix} \H^i(\CC^{\bt}_0(E/F)) & i>0 \\ \H_{-i}\H^0(N\CC^{\bt}_{\bt}(E/F)) & i \le 0.  \end{matrix} \right.
$$
\end{definition}

\begin{lemma}\label{cohochar}
If  $V \in sc\Vect_k$ is quasi-smooth, then the  inclusion map 
$$
\lrcorner V \to \Tot NV 
$$
is a quasi-isomorphism, and 
$$
\H_i(NZ^nV) \cong \H^{n-i}(\lrcorner V)
$$
for all $i, n \ge 0$.
\end{lemma}
\begin{proof}
Combine the proofs of \cite{ddt1} Lemma \ref{ddt1-cohochar} and \cite{ddt1} Proposition \ref{ddt1-totcoho}.
\end{proof}

\begin{lemma}\label{ddcqschar}
A levelwise smooth morphism $f:E \to F$ of pre-DDCs is quasi-smooth  if and only if  $\CC^{\bt}_{\bt}(E/F)$ is quasi-smooth (in the sense of Definition \ref{q12vect}).
\end{lemma}
\begin{proof}
Since $f$ is levelwise smooth, we know by Proposition \ref{smoothchar} that each $E^i \to F^i$ is a smooth map of functors $\C_{\L} \to \bS$. For a small extension $A \to B$ in $\C_{\L}$ with kernel $I$, we thus deduce that 
$E^i(A) \to F^i(A)\by_{F^i(B)}E^i(B)$ is a fibration, with fibre $\CC^i_{\bt}(E/F) \ten I$.

Hence  $\H_*\CC^i_{\bt}(E/F)=0$ for all $i$ if and only if we have $E^i \to F^i$ trivially smooth for all $i$, i.e. if Definition \ref{q12}.(Q1) holds. The result now follows from the characterisation of Definition \ref{q12vect}. 
\end{proof}

\begin{definition}
A  morphism $f:E \to F$ of DDCs is said to be a quasi-isomorphism if $\H^*(\lrcorner f):\H^*(\lrcorner E) \to \H^*(\lrcorner F)$ is an isomorphism.
\end{definition}

\begin{definition}\label{mcdef}
Recall from \cite{ddt2} Definition \ref{ddt2-mcdef} that for any  pre-SDC $E$, we define the  Maurer-Cartan functor  $\mc_E:s\C_{\L}\to \Set$ by
$$
\mc_E(A)\subset \prod_{n\ge 0} E^{n+1}(A^{I^n}),
$$
consisting of those $\underline{\omega}$ satisfying:
\begin{eqnarray*}
\omega_m(s_1,\ldots, s_m)*\omega_n(t_1,\ldots, t_n)&=&\omega_{m+n+1}(s_1,\ldots, s_m,0,t_1,\ldots, t_n);\\ 
\pd^i\omega_n(t_1,\ldots,t_n)&=&\omega_{n+1}(t_1, \ldots,t_{i-1},1,t_i,\ldots,t_n);\\ 
\sigma^i\omega_n(t_1,\ldots,t_n)&=&\omega_{n-1}(t_1, \ldots,t_{i-1},\min\{t_i,t_{i+1}\},t_{i+2},\ldots, t_n);\\
\sigma^0\omega_n(t_1,\ldots,t_n)&=&\omega_{n-1}(t_2, \ldots,t_n);\\
\sigma^{n-1}\omega_n(t_1,\ldots,t_n)&=&\omega_{n-1}(t_1, \ldots,t_{n-1}),\\ 
\sigma^0\omega_0&=&1,
\end{eqnarray*}
where $I:=\Delta^1$. 
\end{definition}

\begin{definition}\label{dmcdef}
Given a pre-DDC $E$, define the derived Maurer-Cartan functor 
$
\MC(E):s\C_{\L} \to \bS
$
by $\Hom_{\bS}(K,\MC(E)):= \mc(E_K)$.
\end{definition}

\begin{proposition}
If $f: E \to F$ is a quasi-smooth  morphism of pre-DDCs, then 
$$
\MC(f):\MC(E) \to \MC(F)
$$
is quasi-smooth, with cohomology groups
$$
\H^i(\MC(E)/ \MC(F))\cong \H^{i+1}(\lrcorner E/F).
$$

In particular, if $E$ is a DDC, then $\MC(E)$ is quasi-smooth.
\end{proposition}
\begin{proof}
By construction, the simplicial  matching maps are given by
$$
\mc(E_n) \to \mc(F_n)\by_{\mc(F_{\pd \Delta^n})}\mc(E_{\pd \Delta^n}) = \mc( E_{\pd\Delta^n}\by_{F_{\pd\Delta^n}}F_n).
$$
Condition (Q1) from  Definition \ref{q12} for $f$ implies that 
\[
E_n \to E_{\pd\Delta^n}\by_{F_{\pd\Delta^n}}F_n
\]
 is a levelwise smooth map of SDCs, so \cite{ddt2} Proposition \ref{ddt2-mcq} implies that $\MC(f)$
satisfies condition (S1) from Definition \ref{scspqsdef}.

We now need to check that the quasi-smooth partial matching maps
$$
\alpha:\mc(E_n)\to \mc(E_{\L^n_k})\by_{\mc(F_{\L^n_k})} \mc(F_n)
$$ 
are smooth. To do this, we verify the criterion of Corollary \ref{cohosmoothchar}. 

Taking the relative version of  \cite{ddt2}  Proposition \ref{ddt2-mch}, we see that  $\H^*(\alpha)$ is cohomology of the cochain complex
$$
\ker(N_c(E_n/F_n)[1] \to N_c(E_{\L^n_k}/F_{\L^n_k});
$$
this is isomorphic to $ N^s_n N_c(E/F)[1]$, which
gives isomorphisms 
$$
 \H^0(\alpha) \cong N_n\z^1\CC(E/F), \quad \H^i(\alpha)\cong N_n\H^{i+1}(E/F) \text{ for } i>0.
$$
Since $N_n\H^i(E/F)=0$  all $i>0$ by condition (Q2'), we see that $\H^i(\alpha)=0$ for all $i>0$. This implies quasi-smoothness of $\MC(f)$.

Now for $i>0$, the calculations above combine with Lemma \ref{cohochar} to give
$$
\H^i(\MC(E)/ \MC(F))= \H^i(\mc(E_0)/ \mc(F_0))=\H^{i+1}(E_0/F_0)=\H^{i+1}(\lrcorner E/F).
$$
For $i\le 0$, 
$$
\H^i(\MC(E)/ \MC(F))= \H_{-i}\H^0(\CC(E/F)[1])=\H_{-i}\z^1(\CC(E/F))= \H^{i+1}(\lrcorner E/F),
$$
since $\CC(E/F)$ is quasi-smooth (in the sense of Definition \ref{q12vect}) 
\end{proof}

\begin{corollary}
If $f:E \to F$ is a quasi-isomorphism of DDCs, then $\MC(f):\MC(E) \to \MC(F)$ is a weak equivalence.
\end{corollary}

\begin{definition}\label{ddef}
Given a pre-DDC $E$, 
note that
$E^0$ acts on $\MC(E)$ by conjugation. Define the derived deformation functor
$
\Def(E):s\C_{\L} \to \bS
$
by
$$
\Def(E):=[\MC(E)/E^0],
$$
the homotopy quotient (as in Definition \ref{quotdefn}).
\end{definition}

\begin{lemma}\label{quotientsfib}
Given  a  simplicial group $G$,  and a  fibration $X \to Y$ of simplicial $G$-sets, if each $X_n$ and $Y_n$  is a free $G_n$-set, then $X/G \to Y/G$ is a fibration in $\bS$.
\end{lemma}
\begin{proof}
Combine \cite{sht} Corollary V.2.7 and \cite{sht} Lemma V.3.7.
\end{proof}

\begin{corollary}\label{Defcts}
If $f: E \to F$ is a quasi-smooth  morphism of pre-DDCs, then 
$$
(f,q):\Def(E) \to \Def(F)\by_{\bar{W}F^0}\bar{W}E^0,
$$
is quasi-smooth, where $\bar{W}G:= G\backslash WG$ is a model for the classifying space  $BG$ of $G$ (as in \cite{sht} \S V.4).

Thus
$$
\Def(f):\Def(E) \to \Def(F) 
$$
is quasi-smooth, and
$$
\MC(E) \to \Def(E)\by_{\Def(F)}\MC(F)
$$
is a weak equivalence.

In particular, if $E$ is a DDC, then $\MC(E) \to \Def(E)$ is a weak equivalence.
\end{corollary}
\begin{proof}
First observe that $E^0 \to F^0$ is trivially smooth, so $\MC(E)\by WE^0 \to \MC(F)\by WE^0$ is quasi-smooth.

Given a surjection $A \to B$ in $s\C_{\L}$, apply Lemma \ref{quotientsfib}, taking 
$$
X=WE^0(A)\by\MC(E)(A), \quad Y= WE^0(A)\by\MC(E)(B)\by_{\MC(F)(B)}\MC(F)(A),
$$
 and $G=E^0(A)$. This shows that $(f,q)$ satisfies condition (S1) from Definition \ref{scspqsdef}. Condition (Q2) follows similarly, so $(f,q)$ is quasi-smooth.

That $\Def(F)$ is quasi-smooth follows from the observation that $\bar{W}E^0 \to \bar{W}F^0$ is  trivially smooth. 

For the final statements, note that $\MC(E) = \Def(E)\by_{\bar{W}E^0}1$, and  
$$
Y:=\Def(E)\by_{\Def(F)}\MC(F)=\Def(E)\by_{\bar{W}F^0}1
$$
If $Z:= \bar{W}E^0\by_{\bar{W}F^0}1$, then $Z$ is trivially fibrant, so $1 \to Z$ is a weak equivalence. The map $\MC(E)\to Y$ is then just the pullback of $1 \to Z$ along $Y \to Z$.
\end{proof}

\subsection{Comparison with SDCs}

\begin{definition}\label{powersdc}
Given a pre-DDC $E$, define a pre-DDC $DE$ by $(DE)_n:= (E_n)^{\Delta^n}$, in the notation of \cite{ddt2} Definition \ref{ddt2-powersdc}, i.e.
$$
(E^X)^n=(E^n)^{X_n}.
$$
For $x \in X_{n+1}$, $y \in Y_{n+1}$, $z \in X_{m+n}$, $1 \le i \le n$, $0 \le j < n$, $e \in (E^X)^n$ and $f \in (E^X)^m$, we define the operations by
\begin{eqnarray*}
\pd^i(e)(x)&:=& \pd^i(e(\pd_i x))\\
\sigma^j(e)(y)&:=& \sigma^j(e(\sigma_i y)),\\ 
(f*e)(z)&:=& f((\pd_{m+1})^nz)*e((\pd_{0})^mz).
\end{eqnarray*}
\end{definition}

\begin{proposition}\label{Dworks}
If $f:E \to F$ is a  map of pre-DDCs with 
\begin{enumerate}
\item $f^i:E^i \to F^i$ smooth for all $i$, and
\item  $\H^i(E/F)$  a constant simplicial complex for all $i>0$,
\end{enumerate}
  then $Df:DE \to DF$ is quasi-smooth. 

In particular, $DE$ is a DDC for all SDCs $E$.
\end{proposition}
\begin{proof}
By smooth base change, we know that $Df$ is levelwise smooth. We now verify the conditions of Lemma \ref{ddcqschar}.
$$
\CC^i(DE/DF)_n= \CC^i(E_n/F_n)^{\Delta^n_i}= \CC^i(E_n/F_n)\ten \CC^i(\Delta^n, k).
$$
Thus
$$
\H_*(DE/DF)^i= \H_*(E^i/F^i)\ten \H_*(\CC^i(\Delta^{\bt}, k))=0,
$$
since 
the cosimplicial complex $k\ten \Delta^{\bt}_i$ is contractible. Moreover,
$$
\H^*(DE/DF)_n= \H^*(E_n/F_n)\ten \H^*(\Delta^n, k)= \H^*(E_n/F_n),
$$
so $\H^i(DE/DF)_n$ is constant for $i>0$.
\end{proof}

\begin{lemma}\label{levelwiseqscoho}
If   $X$ in $sc\Sp$ is a levelwise quasi-smooth object for which the simplicial vector spaces $\H^i(\tan X)$ have constant simplicial structure for $i>0$, then 
$$
\H^*(X) \cong \H^*(\Tot N\tan X),
$$  
as defined in Definitions \ref{generalcoho} and \ref{tandef} respectively.
\end{lemma}
\begin{proof}
By \cite{ddt1} Lemma \ref{ddt1-levelwiseqs}, there is a weak equivalence $X \to \underline{X}$ to a quasi-smooth object. Since this is also a levelwise weak equivalence, we find that $\Tot N\tan X \to \Tot N\tan\underline{X}$ is a quasi-isomorphism. Now Proposition \ref{ctodtnew} implies that $\underline{X}$ is quasi-smooth, so
 $\H^*(X) \cong \H^*(\underline{X}) \cong \H^*(\Tot N\tan X)$, the last isomorphism coming from \cite{ddt1}  Theorem \ref{ddt1-totcoho}.
\end{proof}

\begin{proposition}\label{Dworkswell}
If $E$ is a pre-DDC for which $E \to \bt$  satisfies the conditions of Proposition \ref{Dworks}, then $\alpha:\Def(E)\to\Def(DE)$ is a quasi-smooth replacement for $\Def(E)$.
\end{proposition}
\begin{proof}
By Proposition \ref{Dworks}, we know that $\Def(DE)$ is quasi-smooth, so we just need to show that $\alpha$ is a weak equivalence. Now, for $i>0$,
$$
\H^i(\tan \Def(E))= H^i(\tan \MC(E))= \H^{i+1}(E),
$$
which has constant simplicial structure, so $\H^i(\Def(E))= \H^i(\Tot N\tan\Def(E))= \H^{i+1}(\Tot N\CC^{\bt}(E))$.

Similarly, $\H^i(\Def(DE))= \H^{i+1}(\Tot N\CC^{\bt}(DE))\cong \H^{i+1}(\Tot N\CC^{\bt}(E))$, so  \cite{ddt1} Corollary \ref{ddt1-weakchar} ensures that $\alpha$ is a weak equivalence.
\end{proof}

\begin{corollary} 
For an SDC $E$, the functors  $\ddef_E$ (from \cite{ddt2} Definition \ref{ddt2-ddefdef})
 and $\Def(DE)$ (and hence $\MC(DE)$) are weakly equivalent (in $sc\Sp$).
\end{corollary}
\begin{proof}
Recall that $\ddef_E$ is $[\mmc_E/\uline{E}^0]$, which is a 
 quasi-smooth replacement  of $[\mc_E/E^0]$ by \cite{ddt1} Lemma \ref{ddt1-levelwiseqs}, in the sense that there is a weak equivalence $[\mc_E/E^0] \to \ddef_E$. If we let $E$ denote the constant pre-DDC $E_n:=E$, then $[\mc_E/E^0] =\Def(E)$, and we may apply Proposition \ref{Dworkswell}.
\end{proof}

\begin{lemma}\label{pioD}
If $E$ is a pre-DDC for which $E \to \bt$  satisfies the conditions of Proposition \ref{Dworks}, then for all $ A \in \C_{\L}$,  the map $\alpha(A):\Def(E)(A)\to\Def(DE)(A)$ is a weak equivalence in $\bS$. 
\end{lemma}
\begin{proof}
First observe that $\Def(E)(A)_0 =\mc_{E_0}(A)=\Def(DE)(A)$, and write $\pi^0F:= F|_{\C_{\L}}$. Now, 
$
\tan \pi^0 \Def(E)
$
is the mapping cone of $\CC^0(E) \xra{d_c} \z^1(E)$, so (Q2) ensures that $\pi_n \tan \pi^0\Def(E)=\H^{1-n}(\lrcorner E)$. Observe that $\pi_n \tan \pi^0\Def(DE)=\H^{1-n}(\lrcorner E)$, similarly.

The proofs of  Proposition \ref{robs} and Corollary \ref{weak} adapt to show that $\alpha(A)$ is a weak equivalence in $\bS$ for all $A$, by taking  small extensions $A \to B$ with kernel $I$, and considering the long exact  sequences
$$
\begin{CD}
\ldots @>>> \H^{1-n}(\lrcorner E)\ten I @>>>\pi_n(  \Def(E)(A),x) @>>> \pi_n(\Def(E)(B),x) @>>>  \ldots\\
@VVV @VVV @VVV @VVV @VVV\\
\ldots @>>> \H^{1-n}(\lrcorner E)\ten I @>>> \pi_n(\Def(DE)(A),x) @>>> \pi_n(\Def(DE)(B),x) @>>>  \ldots .
\end{CD}
$$
associated to the fibrations $\Def(E)(A)\to \Def(E)(B)$ and $\Def(DE)(A) \to \Def(DE)(B)$.
\end{proof}

\section{Constructing DDCs}\label{constrsn}

\subsection{Simplicial monadic adjunctions}

A simplicial category $\C$ has a class $\Ob\C$ of objects, and for all $A,B \in \Ob\C$, a simplicial set $\underline{\Hom}_{\C}(A,B)$ of morphisms, with the usual multiplication and identity properties. 
For a simplicial category $\C$, we denote by $\C_n$ the category with objects $\Ob\C$ and morphisms $\Hom_{\C_n}(A,B):= \underline{\Hom}_{\C}(A,B)_n$. 

\begin{definition}
Say that a functor $F:\C \to \cD$ of simplicial categories is an equivalence if the functors $F_n:\C_n \to \cD_n$ are all equivalences.
\end{definition}

\begin{definition}
Given a simplicial category $\C$, set $\Hom_{\C}(A,B):= (\underline{\Hom}_{\C})_0(A,B)$. 
\end{definition}

\begin{definition}
For simplicial categories $\cD, \cE$, and a pair of functors
$$
\xymatrix@1{\cD \ar@<1ex>[r]^G & \cE \ar@<1ex>[l]^F},
$$
recall that an adjunction $F \dashv G$ is a natural isomorphism
$$
\underline{\Hom}_{\cD} (FA,B) \cong \underline{\Hom}_{\cE}(A,GB).
$$
We say that $F$ is left adjoint to $G$, or $G$ is right adjoint to $F$. Let $\bot = FG$, and $\top=GF$. To give an adjunction is equivalent to giving two natural transformations, the unit and co-unit
$$
\eta:\id_{\cE} \to  \top,\quad \vareps: \bot \to \id_{\cD},
$$
satisfying the triangle identities $\vareps F \circ F\eta =\id_F $, $G\vareps \circ\eta G =\id_G$.  
\end{definition}

Given an adjunction
$$
\xymatrix@1{\cD \ar@<1ex>[r]^U_{\top} & \cE \ar@<1ex>[l]^F}
$$
with unit $\eta:\id \to UF$ and co-unit $\vareps:FU \to \id$, we let $\top=UF$, and define the simplicial category  $\cE^{\top}$ of $\top$-algebras to have objects
$$
\top E \xra{\theta} E,
$$
for $\theta \in \underline{\Hom}_0(\top E,E)$, such that $\theta\circ \eta_E=\id$ and $\theta \circ \top \theta= \theta \circ  U\vareps_{FE}$.
We define morphisms by setting
$$
\underline{\Hom}_{\cE^{\top}}(\top E_1 \xra{\theta} E_1, \top E_2 \xra{\phi} E_2 )\subset \underline{\Hom}_{\cE}(E_1,E_2)
$$
to be the equaliser of
$$
\xymatrix@1{ \underline{\Hom}_{\cE}(E_1,E_2) \ar@<0.5ex>[r]^{\phi_*\circ \top}  
\ar@<-0.5ex>[r]_{\theta^*} 
& \underline{\Hom}_{\cE}(\top E_1,E_2).}
$$

We define the comparison functor $K:\cD \to \cE^{\top}$ by
$$
B \mapsto ( UFUB \xra{U\vareps_B}  UB )
$$
on objects, and $K(g)=U(g)$ on morphisms.

\begin{definition}
An adjunction 
$$
\xymatrix@1{\cD \ar@<1ex>[r]^U_{\top} & \cE \ar@<1ex>[l]^F},
$$
of simplicial categories is said to be  monadic if $K:\cD \to \cE^{\top}$ is an equivalence.
\end{definition}

\begin{examples}
Intuitively, monadic adjunctions correspond to algebraic theories, such as the adjunction
$$
\xymatrix@1{\Ring \ar@<1ex>[r]^U_{\top} &\Set  \ar@<1ex>[l]^{\Z[-]},}
$$
 between rings and sets, $U$ being the forgetful functor. Other examples are $k$-algebras   over $k$-vector spaces, or groups over sets.
\end{examples}

\begin{definition}
Given an adjunction
$$
\xymatrix@1{\cD \ar@<1ex>[r]^V & \cE \ar@<1ex>[l]^G_{\bot}},
$$
let $\bot=VG$, so $\bot^{\op}$ is a monad on $\cE^{\op}$. Define $\cE_{\bot}:= ((\cE^{\op})^{\bot^{\op}})^{\op}$, with $K=K^{\op}: \cD \to \cE_{\bot}$.
The adjunction  is said to be  comonadic  if $K:\cD \to \cE_{\bot}$ is an equivalence.
\end{definition}

\begin{example}
If $X$ is a topological space (or any site with enough points) and $X'$ is the set of points of $X$, let  $u:X'\to X$ be the associated morphism.
Then the adjunction $u^{-1}\dashv u_*$ on sheaves is comonadic, so the category of sheaves  on $X$ is equivalent  $u^{-1}u_*$-coalgebras in  the  category of sheaves  (or equivalently presheaves)  on $X'$

A more prosaic example is that for any ring $A$, the category of $A$-coalgebras is comonadic over the category of $A$-modules.
\end{example}

\subsubsection{Bialgebras }\label{bialgsn}
As in \cite{osdol} \S IV, take a category $\cB$ equipped with both a monad $(\top, \mu, \eta)$ and a comonad $(\bot, \Delta, \gamma)$, together with a distributivity transformation
$\lambda: \top\bot \abuts \bot \top$ satisfying various additional conditions.

\begin{definition}\label{bialgdef}
Given a distributive monad-comonad pair $(\top, \bot)$ on a simplicial category $\cB$, define the category $\cB^{\top}_{\bot}$ of bialgebras as follows. The objects of  $\cB^{\top}_{\bot}$ are triples $(\theta, B, \beta)$ with $(\top B \xra{\theta} B)$ an object of $\cB^{\top}$ and $B \xra{\beta} \bot B$ an object of $\cB_{\bot}$, such that the composition  $(\beta \circ \theta):\top B \to \bot B$ agrees with the composition 
$$
\top B \xra{\top \beta} \top\bot B \xra{\lambda} \bot \top B \xra{\bot \theta} \bot B. 
$$  
Morphisms are then given by setting
$$
\underline{\Hom}_{\cB^{\top}_{\bot}}(\top B \xra{\theta} B \xra{\beta} \bot B, \top B' \xra{\theta'} B' \xra{\beta'} \bot B' )\subset \underline{\Hom}_{\cB}(B,B')
$$
to be the equaliser of
$$
\xymatrix@1{ \underline{\Hom}_{\cB}(B,B') \ar@<0.5ex>[rr]^-{(\theta'_*\circ \top, \beta'_*)}  
\ar@<-0.5ex>[rr]_-{(\theta^*, \beta^*\circ \bot)} 
&& \underline{\Hom}_{\cB}(\top B,B')\by \underline{\Hom}_{\cB}(B, \bot B').}
$$
\end{definition}

\begin{example}
If $X$ is a topological space (or any site with enough points) and $X'$ is the set of points of $X$, let $\cD$ be the category of sheaves  of rings on $X$. If $\cB$ is the category of sheaves (or equivalently presheaves) of sets on $X'$, then the description above characterises $\cD$ as a category of bialgebras over $\cB$, with the comonad being $u^{-1}u_*$ for $u:X'\to X$, and the monad  being the free polynomial  functor. 
\end{example}

\subsection{The construction}\label{construct}

We let $s\Cat$ denote the category of  simplicial categories. 

\begin{definition}
Given functors $\cA \xra{f} \cB \xla{g} \C$ of simplicial categories, define the fibre product $\cA\by_{\cB}\C$  by
$$
\Ob(\cA\by_{\cB}\C)=\{(A,\beta,C)\,:\, A \in \Ob \cA,\, C \in \Ob C, \,\beta \in \Iso_{\cB_0}(fA,gC)\},
$$
with morphisms 
$$
\underline{\Hom}_{\cA\by_{\cB}\C}((A,\beta,C), (A',\beta',C'))= \underline{\Hom}_{\cA}(A,A')\by_{\beta'_*f, \underline{\Hom}_{\cB}(fA,gC'), \beta^*g}\underline{\Hom}_{\C}(C,C').
$$
\end{definition}

\begin{definition}
We say that a functor $F: \C_{\L}\to \Set$ is homogeneous if for all small extensions $A \to B$ in $\C_{\L}$, 
$$
F(A\by_BC)\to F(A)\by_{F(B)}F(C)  
$$
is an isomorphism. Note that this is equivalent to being a disjoint union of left-exact functors.

Similarly, a functor $\cD: \C_{\L}\to s\Cat$ is said to be homogeneous if  
$$
F(A\by_BC)\to F(A)\by_{F(B)}F(C)  
$$
is an equivalence for all small extensions $A \to B$.
\end{definition}

\begin{definition}\label{uniftriv}
We say that a homogeneous functor $\cB:\C_{\L}\to s\Cat$   has uniformly trivial deformation theory  if  
\begin{enumerate}
\item for all $A \in \C_{\L}$ and all $B_1,B_2 \in \Ob \cB(A)$, the functor $\underline{\Hom}_{\cB}(B_1,B_2):\C_A \to \Set$ of morphisms from $B_1$ to $B_2$ is trivially smooth (in the sense of Definition \ref{pioqs}) and  homogeneous;
\item  for $A' \onto A$ in $\C_{\L}$, $\cB_0(A') \to \cB_0(A)$ is essentially surjective. 
\end{enumerate}
\end{definition}

Now, assume that we have a diagram
$$
\xymatrix@=8ex{
\cD \ar@<1ex>[r]^{U}_{\top} \ar@<-1ex>[d]_{V} 
&\ar@<1ex>[l]^{F} \cE  \ar@<-1ex>[d]_{V} 
\\
\ar@<-1ex>[u]_{G}^{\dashv}	\cA \ar@<1ex>[r]^{U}_{\top} 
&\ar@<1ex>[l]^{F} \ar@<-1ex>[u]_{G}^{\dashv} \cB,  
}
$$
of adjunctions of homogeneous simplicial category-valued functors on $\C_{\L}$, 
with $F\dashv U$ monadic and  $G\vdash V$  comonadic. Let 
\begin{align*}
\toph&=UF&	\both&=FU\\
\botv&=VG&	\topv&=GV,
\end{align*}
with 
$$
\eta:1 \to \toph,\quad \gamma:\botv \to 1,\quad \vareps:\both \to 1 \text{ and }\alpha:1 \to \topv.
$$
Assume that these adjunctions satisfy the simplicial analogues of \cite{osdol} \S IV or \cite{paper2} \S \ref{paper2-gensdc}, in other words that $U$ and $V$ commute with everything (although $G$ and $F$ need not commute). 

Fix $D \in \Ob \cD(k)$, such that we may lift $UVD \in \Ob\cB(k)$ to $B \in \Ob\cB(\L)$, up to isomorphism (in $\cB_0(k)$).

\begin{theorem}\label{ddcmain}
There is a natural pre-DDC $E$ associated to this diagram, given by
$$
E^n=\underline{\Hom}_{\cB}(\toph^n B, \botv^n B)_{UV(\alpha^n_{D}\circ \vareps^n_D)}.
$$

If $E$ is 
levelwise smooth, satisfying Condition (Q2) of Definition \ref{q12},
then  the classifying space $\bar{W}\cD_{D,\id}$  is canonically weakly equivalent to the restriction  $\pi^0\Def(E)$ (from Lemma \ref{pioD}) as a functor from  $\C_{\L}$ to $\bS$. Here $\cD_{D,\id}(A)$ is the simplicial groupoid given by the fibre product 
$$
\cD(A)\by_{\cD(k)}(D,\id),
$$
where $(D,\id)$ is the category with one object and one morphism. 
\end{theorem}
\begin{proof}
For each $m$, $E_m$ is the SDC defined in \cite{paper2} \S \ref{paper2-gensdc} associated to the monad $\toph$ and comonad $\botv$ over the category $\cB_m$.

Since the adjunctions are monadic or comonadic, the proof of \cite{paper2} Theorem \ref{paper2-main} adapts to give functorial equivalences
$$
K(A): \cD(A) \to \cB^{\toph}_{\botv}(A)
$$
between $\cD$ and the simplicial category of $(\toph,\botv)$-bialgebras.

Let $D' \in \Ob \cB^{\toph}_{\botv}(k)$ be the bialgebra over $\bar{B} \in \Ob \cB(k)$, with  bialgebraic structure coming from the isomorphism $UVD \cong \bar{B}$. Let $\cG$ be the full subcategory of $\cB^{\toph}_{\botv}(A)$ on objects
$$
\mc_{E_0}(A)= \{ \omega \in \Ob \cB^{\toph, \botv}(A)\,:\, \bar{X}=D' \in \Ob \cB^{\toph, \botv}(k).
$$
Morphisms in $\cG$ are just
$$
\underline{\Hom}_{\cG}(\omega, \omega')= \{f \in E^0\,:\, f*\omega=\omega'*f\},
$$
from which we deduce that $\cG$ is a simplicial groupoid. Moreover, observe that $\cG \to(\cB^{\toph}_{\botv})_{D', \id}$ is an isomorphism of simplicial categories, so $\cG$ is equivalent to $\cD_{D,\id}$. In particular, this implies that $\cD_{D,\id}$ is a simplicial groupoid. It therefore suffices to compare $\cG$ with $\Def(E)$.

\begin{lemma}
The functor $\cG:\C_{\L} \to s\gpd$ is quasi-smooth, in the sense that it maps small extensions to fibrations (as defined in \cite{sht} \S V.7).
\end{lemma}
\begin{proof}[Proof of lemma]
Smoothness of $E^0_0$ implies that the path-lifting property is satisfied. Given $K \into L \in \bS$, and a small extension $A \to B$ with kernel $I$, the obstruction to lifting the diagram
$$
\begin{CD}
K @>>> \underline{\Hom}_{\cG}(\omega, \omega')(A)\\
@VVV @VVV\\
L @>>> \underline{\Hom}_{\cG}(\omega, \omega')(A)
\end{CD}
$$
lies in $\H^1( \ker(\CC^{\bt}(E_L) \to \CC^{\bt}(E_K)))$. If we write $V^{\bt}= \ker(\CC^{\bt}(E_L) \to \CC^{\bt}(E_K))$, then we have an exact sequence
$$
\H^0(E_L) \xra{\alpha} \H^0(E_K)\to \H^1(V^{\bt}) \to \H^1(E_L) \xra{\beta} \H^1(E_K).
$$
If $K \into L$ is a trivial cofibration, then Condition (Q2) of Definition \ref{q12} ensures that $\alpha$ is surjective and $\beta$ an isomorphism, so the obstruction is zero and the lift exists, proving that $\underline{\Hom}_{\cG}(\omega, \omega')$ is quasi-smooth, as required.  
\end{proof}

Now, the inclusions $\mc_{E_0} \into \mc(E)$ and $\underline{\Hom}_{\cG}(\omega, \omega') \into E^0$ define a morphism
$$
\alpha:\bar{W}\cG \to \pi^0\Def(E)
$$
of quasi-smooth functors $\C_{\L} \to \Set$. Note that $\bar{W}\cG_0 = \mc_{E_0}= \pi^0\Def(E)_0$.

(Q2) also ensures that $\alpha$ is a weak equivalence on tangent spaces, with $\pi_n \bar{W}\cG(k[\eps])= \H^{1-n}(\lrcorner E)$. As in the proof of Lemma \ref{pioD}, this implies  that $\alpha(A)$ is a weak equivalence in $\bS$ for all $A$.
\end{proof}

\begin{remark}
If  $\cB$ has uniformly trivial deformation theory, then 
note that $\tilde{D}$ always exists, and that the pre-DDC $E$ of Theorem \ref{ddcmain} automatically satisfies Definition \ref{q12}.(Q1).

However, if $E$ just satisfies all the conditions of  Theorem \ref{ddcmain}, then  
Propositions \ref{Dworks} and \ref{Dworkswell} then give a DDC $DE$, which by Lemma \ref{pioD} also has $\bar{W}\cD_{D,\id} \sim \pi^0\Def(DE)$.
\end{remark}

\subsection{Deformations of diagrams and invariance under weak equivalence}\label{weakinvart}

In a similar vein, we may study deformations of a morphism, or even of a diagram. 

\begin{definition}\label{delta**}
Define $\Delta_{**}$ to be the subcategory of the ordinal number category $\Delta$ containing only  those morphisms $f:\mathbf{m} \to \mathbf{n}$ with $f(0)=0, f(m)=n$.  Given a category $\C$, a functor $X:\Delta_{**} \to \C$ consists of objects $X^n \in \C$, with all of the operations $\pd^i, \sigma^i$ of a cosimplicial complex except $\pd^0, \pd^{n+1}:X^n\to X^{n+1}$. 
\end{definition}

\begin{definition}
Given a monoidal category $\C$ and a set $\cO$, recall from \cite{monad} that a $\C$-valued quasi-descent datum $X$ on objects $\cO$ consists of:
\begin{enumerate}
		\item objects $X(a,b) \in \C^{\Delta_{**}}$ for all $a,b \in \cO$;
		\item morphisms $X(a,b)^m \ten X(b,c)^n \xra{*} X(a,c)^{m+n}$ making the following diagram commute for all $,b,c \in \cO$
		$$
		\begin{CD} \Delta_{**} \by \Delta_{**} @>>{ X(a,b)\ten X(b,c)}> \C  \\
		@V{\by }VV @VV{*}V\\
		\Delta_{**} @>>{X(a,c) }> \C.
		\end{CD}
		$$
		\item morphisms $1 \to X(a,a)^0$    for all $a \in \cO$, acting as the identity for the multiplication $*$.
		\end{enumerate}
\end{definition}

Note that a pre-DDC over $\L$ is a quasi-descent datum (on one object) in the monoidal category $(s\Sp, \by)$. 

\begin{definition}
Let $Q\Dat(\C)$ be the category of $\C$-valued quasi-descent data, i.e. of pairs $(\cO,X)$ for $\cO$ a set and $X$ a quasi-descent datum on objects $\cO$. 

We say that $\cD$ is an enrichment of a $\C$-enriched category $\cF$ if $\Ob \cF \cong \Ob \cD$ and $\cF(x,y) \cong \cD^0(x,y)$, compatible with the product and identities.
\end{definition}

\begin{proposition}\label{enrich}
For a diagram of simplicial category-valued functors as in \S \ref{construct}, the  $s\Sp$-enriched category $\cB(\L)$ has a natural  enrichment in $Q\Dat(s\Sp)$. If the simplicial structure on $\cB$ is constant, then  this enrichment is in $Q\Dat(\Sp)$.
\end{proposition}
\begin{proof}
This is just \cite{monad} Proposition \ref{monad-enrichtopbot}.
The enriched $\Hom$-set $\underline{\hom}(B,B'):\C_{\L} \to \bS^{\Delta_{**}}$ is given by
$$
\underline{\hom}^n(B,B'):=\underline{\Hom}_{\cB}(\toph^n B, \botv^n B').
$$

If the simplicial structure on $\cB$ is constant, then $\underline{\Hom}_{\cB}= \Hom_{\cB}$, so $\underline{\hom}(B,B')$ lies in $\Sp$.
\end{proof}

\begin{definition}
Given a morphism $f:D \to D'$ in $\cD(k)$  for which $UVD, UVD'$ lift to $B,B'$ in $\cB(\L)$, define
$$
E^n_{\cD/\cB}(f):=\underline{\hom}^n(B,B')_{UV(\alpha^n_{D'}\circ f \circ \vareps^n_D)} \in s\Sp
$$

Write $E^*_{\cD/\cB}(D):= E^*_{\cD/\cB}(\id_D)$.
\end{definition}

\begin{definition}
Given a morphism $f:D \to D'$ in $\cD(k)$ for which $E^*_{\cD/\cB}(f) \in (s\Sp_k)^{\Delta_{**}}$ is levelwise smooth, define 
$$
\CC^{\bt}_{\cD/\cB}(f):=\tan E^*_{\cD/\cB}(f),
$$
and note that that this becomes a cosimplicial complex (of simplicial complexes), by \cite{ddt2} Lemma \ref{ddt2-cosimplicialE}. Explicitly,
$$
\CC^n_{\cD/\cB}(f)= \tan \HHom_{\cD|_{\C_k}}(\both^{n+1}D,\toph^{n+1} D')_{\alpha^{n+1}_{D'}\circ f \circ \vareps^n_D}.
$$

Define 
\begin{eqnarray*}
\underline{\Ext}^n_{\cD/\cB}(f)&:=& \H^n(\CC^{\bt}_{\cD/\cB}(f)) \in s\Vect_k\\ 
\Ext^n_{\cD/\cB}(f)&:=& \H^n(\Tot N\CC^{\bt}_{\cD/\cB}(f)) \in \Vect_k.
\end{eqnarray*}
\end{definition}

\begin{definition}\label{qsmorphism}
Say that a morphism $f:D \to D'$ in $\cD(k)$ is Q2 over $\cB$ if 
\begin{enumerate}
\item $UVD, UVD'$ lift to $\cB(\L)$,

\item  $E^*_{\cD/\cB}(f)$ is levelwise smooth, and

\item $\underline{\Ext}^i_{\cD/\cB}(f)$
is a constant simplicial complex for $i>0$. 
\end{enumerate}

We say that $f$ is quasi-smooth over $\cB$ if in addition $\H_*\CC^{n}_{\cD/\cB}(f)=0$ for all $n$.
\end{definition}

\begin{remark}
Note that if $f$ is Q2 over $\cB$,  then we have $\Ext^*_{\cD/\cB}(f)=\H^*(\lrcorner N\CC^{\bt}_{\cD/\cB}(f))$, by Lemma \ref{cohochar}.
\end{remark}

\begin{definition}\label{sdcdiagram}
Given a small category $\bI$, and an $\bI$-diagram $\bD:\bI \to \cD(k)$ with objects $UV\bD(i)$ lifting to  $\cB(\L)$,
define the pre-DDC $E^{\bt}_{\cD/\cB}(\bD)$ by
$$
E^n_{\cD/\cB}(\bD)= \prod_{\substack{i_0 \xra{f_1} i_1 \xra{f_2} \ldots \xra{f_n} i_n\\ \text{in }\bI}} E^n(\bD(f_n\circ f_{n-1}\circ \ldots f_0)) = \prod_{x \in B\bI_n} E^n( \bD(\pd_1^{\phantom{1}n-1}x)),
$$
where $B\bI$ is the nerve of $\bI$ (so $B\bI_0=\Ob(\bI)$,  $B\bI_1=\Mor(\bI)$), and  $\pd_1^{\phantom{1}-1}:=\sigma_0$.
The operations are defined as in Definition \ref{powersdc}.
\end{definition}

\begin{lemma}\label{q2diagram}
Given an $\bI$-diagram $\bD:\bI \to \cD(k)$ with all morphisms $\bD(f)$ quasi-smooth (resp. Q2) over $\cB$, the pre-DDC $E^{\bt}_{\cD/\cB}(\bD)$ is quasi-smooth (resp. is levelwise smooth and satisfies Definition \ref{q12}.(Q2)).
\end{lemma}

\begin{proposition}\label{governdiagram}
Given an $\bI$-diagram $\bD:\bI \to \cD(k)$ with all morphisms Q2,   the classifying space $\bar{W}(\cD^{\bI})_{\bD,\id}$ and  $\pi^0\Def(E^{\bt}_{\cD/\cB}(\bD) )$ are canonically weakly equivalent as functors from  $\C_{\L}$ to $\bS$.
\end{proposition}
\begin{proof}
This is just \cite{monad} Lemma \ref{monad-sdcdiagramsub}.
\end{proof}

\begin{definition}
Say that a morphism $f:D \to D'$ in $\cD(k)$ is an $\Ext_{\cD/\cB}$-equivalence if $UVD, UVD'$ lift to $B,B'$ in $\cB(\L)$, with $E^*(f)$ levelwise smooth, and  the maps
$$
\Ext^*_{\cD/\cB}(\id_D) \xra{f_*} \Ext^*_{\cD/\cB}(f) \xla{f^*} \Ext^*_{\cD/\cB}(\id_{D'})
$$
are isomorphisms.
\end{definition}

\begin{proposition}\label{weakin}
If a morphism $f:D \to D'$ in $\cD(k)$ is an $\Ext_{\cD/\cB}$-equivalence, with the morphisms $f,\id_D, \id_{D'}$ all Q2,   then the DDCs
$
DE^{\bt}_{\cD/\cB}(D)
$
and $
DE^{\bt}_{\cD/\cB}(D')
$
are quasi-isomorphic.
\end{proposition}
\begin{proof}
Let $\bI:= (\bt \to \bt)$ be the category with two objects and one non-identity morphism, and consider the diagram $\bD:\bI \to \cD(k)$ given by  $D \xra{f} D'$.
By Lemma \ref{q2diagram} and Proposition \ref{Dworks}, we know that the pre-DDCs $ DE^{\bt}_{\cD/\cB}(D), DE^{\bt}_{\cD/\cB}(D')$ and $DE^{\bt}_{\cD/\cB}(\bD)$ are all DDCs.

The inclusions of objects into  $\bI$ give morphisms $E^{\bt}_{\cD/\cB}(D) \la E^{\bt}_{\cD/\cB}(\bD)\to E^{\bt}_{\cD/\cB}(D')$. We just need to describe the cohomology groups $\H^*(\lrcorner E^{\bt}_{\cD/\cB}(\bD))$ to show that these induce quasi-isomorphisms.

The tangent space $\CC^{\bt}_{\cD/\cB}(\bD)$ is the diagonal cosimplicial complex associated to the bicosimplicial complex
$$
\prod_{x \in B\bI_m} \CC^n_{\cD/\cB}( \bD(\pd_1^{\phantom{1}m-1}x)),
$$
whose horizontal normalisation is the cochain complex 
$$
\CC^{\bt}_{\cD/\cB}(\id_D) \by \CC^{\bt}_{\cD/\cB}(\id_{D'}) \xra{(f_*,-f^*)}  \CC^{\bt}_{\cD/\cB}(f)
$$
in degrees $0$ and $1$.

Thus $\CC^{\bt}_{\cD/\cB}(\bD)$ is the mapping cone of the morphism $(f_*, -f^*)$.
Since $f$ is an $\Ext_{\cD/\cB}$-equivalence, we deduce that the maps $\Ext^*_{\cD/\cB}(\id_D) \la  \H^*(\lrcorner E^{\bt}_{\cD/\cB}(\bD))\to \Ext^*_{\cD/\cB}(\id_{D'})$ are indeed isomorphisms.
\end{proof}

\subsubsection{Constrained deformations}

We now consider a generalisation, by taking a small diagram 
$$
\bD: \bI\to \cD(k), 
$$ 
 a subcategory  $\bJ \subset \bI$, and $\widetilde{\bD|_{\bJ}}:\bJ \to \cD(\L)$ lifting $\bD|_{\bJ}$. We wish to describe deformations of $\bD$ which agree with $\widetilde{\bD|_{\bJ}}$ on $\bJ$. Note that when $\bI=(0 \to 1)$ and $\bJ=\{1\}$, this is the type of problem considered in  \cite{cones} and \cite{ranlie}.

\begin{proposition}
Given an $\bI$-diagram $\bD:\bI \to \cD(k)$ with all morphisms Q2, and with $\widetilde{\bD|_{\bJ}}$ as above,  the simplicial groupoid of deformations of $\bD$ fixing $\widetilde{\bD|_{\bJ}}$ is governed by the pre-DDC
$$
E^{\bt}_{\cD/\cB}(\bD)\by_{E^{\bt}_{\cD/\cB}(\bD|_{\bJ} )}\bt,
$$
where $
\bt \to E^{\bt}_{\cD/\cB}(\bD|_{\bJ} )
$
is defined by the object of $\mc(E^{\bt}_{\cD/\cB}(\bD|_{\bJ} ))_0 $ corresponding to $\widetilde{\bD|_{\bJ}} $. 
\end{proposition}
\begin{proof}
We need to show that 
the classifying space 
$$
\bar{W}(\cD^{\bI}\by_{\cD^{\bJ}}^h \widetilde{\bD|_{\bJ}})_{\bD,\id}
$$ 
of the homotopy fibre  of simplicial categories is canonically weakly equivalent to 
$$
\pi^0\Def(E^{\bt}_{\cD/\cB}(\bD)\by_{E^{\bt}_{\cD/\cB}(\bD|_{\bJ} )}\bt) 
$$ 
as a functor  from  $\C_{\L}$ to $\bS$.

We know that the functor $\Def$ preserves inverse limits, so 
$$
\pi^0\Def(E^{\bt}(\bD)\by_{E^{\bt}(\bD|_{\bJ} )}\bt)= \pi^0\Def(E^{\bt}_{\cD/\cB}(\bD)) \by_{\pi^0\Def(E^{\bt}(\bD|_{\bJ}))}\bt 
$$

By  Lemma \ref{pioD}, Lemma \ref{q2diagram} and Corollary \ref{Defcts}, we know that $\pi^0\Def(E^{\bt}_{\cD/\cB}(\bD))(A)\to  \pi^0\Def(E^{\bt}(\bD|_{\bJ}))(A)$ is a fibration in $\bS$, so   the fibre over any   point is the homotopy fibre.  Proposition \ref{governdiagram} now shows that this is equivalent to the homotopy fibre of $\bar{W}(\cD(A)^{\bI})\to \bar{W}(\cD(A)^{\bJ})$ over $\widetilde{\bD|_{\bJ}}$, as required.
\end{proof}

\section{Examples}\label{egsn}
We now show how to apply Theorem \ref{ddcmain}, combining it with Definitions \ref{dmcdef} or \ref{ddef} to obtain derived deformation functors. This gives many new examples coming from categories with non-trivial simplicial structure.

\subsection{Chain complexes}\label{chainsn}

We will denote chain complexes by $V_{\bt}$, and their underlying graded modules by $V_*$.

\begin{definition} Define $dg\FMod(A)$ to be  the category  of  chain complexes of flat modules over $A$. We make this into a simplicial category by defining the simplicial normalisation $N^s\underline{\Hom}(U_{\bt},V_{\bt})$ to be the chain complex
$$
N^s_n\underline{\Hom}(U_{\bt},V_{\bt}):=\left\{\begin{matrix} \Hom(U_{\bt},V_{\bt})  & n=0 \\  \prod_{i\ge 0} \Hom(U_i,V_{i+n}) & n>0 \end{matrix},     \right.
$$
with boundary map $\ad_d(f):= d\circ f\pm f\circ d$. This determines the simplicial module $\underline{\Hom}(U_{\bt},V_{\bt}):=(N^s)^{-1}N^s\underline{\Hom}(U_{\bt},V_{\bt})$ by the Dold-Kan correspondence.
\end{definition}

\begin{definition}
Define $g\FMod(A)$ to be  the category  of   flat $\N_0$-graded modules over $A$, with the simplicial 
 structure
$$
N^s_n\underline{\Hom}(U_*,V_*):=\left\{\begin{matrix} \Hom(U_*,V_*)  & n=0 \\   \Hom(U_*,V_*[n])\by\Hom(U_*, V_*[n-1]) & n>0 \end{matrix},     \right.
$$
where the boundary map is given by $d(f,g):=(g,0)$.
\end{definition}

\begin{lemma}
The functor  $g\FMod:\C_{\L}\to s\Cat$ has uniformly trivial deformation theory.  
\end{lemma}
\begin{proof}
Since flat $A$-modules are free, it follows that objects lift. The other properties from Definition \ref{uniftriv} now follow by a simple calculation. 
\end{proof}

\begin{definition}
Let the forgetful functor $dg\FMod \to g\FMod$ be given  by $V_{\bt} \mapsto V_*$, and defined on simplicial morphisms by mapping $f \in N^s_n\underline{\Hom}(U_{\bt},V_{\bt})$ to $(f, \ad_d(f))\in N^s_n\underline{\Hom}(U_*,V_*)$. 
\end{definition}

\begin{lemma}
The forgetful functor $dg\FMod \to g\FMod$ of simplicial categories has a right adjoint $G$, 
and the resulting adjunction is comonadic.
\end{lemma}
\begin{proof} Define
$
(G V_*)_n := V_n \oplus V_{n-1},
$
with $d(v,w)=(w,0)$.  The unit $\alpha:U_{\bt} \to G(U_*)$ of the adjunction is $\alpha(u)=(u,du)$, for any chain complex $U_{\bt}$, and the co-unit $\gamma:G(V_*) \to V_*$ is the map $\gamma(v,w)=v$ of graded modules.
\end{proof}

Let $\bot=VG$ and $\top=GV$, for $V$ the forgetful functor.

\begin{proposition}\label{dgdeform}
For $U_{\bt} \in dg\FMod(k)$,  the pre-DDC 
$$
E^n(A):= \underline{\Hom}_{g\FMod_A}(  \tilde{U}_*\ten A, \bot^n \tilde{U}_*\ten A)_{V(\alpha^n_{U})}
$$
of Theorem \ref{ddcmain} is quasi-smooth, with cohomology
\begin{eqnarray*}
\H^*(\lrcorner E) &=&\H^*( \ldots \xra{\ad_d} \Hom_{g_{\Z}\Mod_k}( U_*, U_*[-n]) \xra{\ad_d} \Hom_{g_{\Z}\Mod_k}( U_*, U_*[-n-1])\xra{\ad_d} \ldots)\\
&=&\EExt^*_{dg_{\Z}\Vect_k}(U_{\bt},U_{\bt}),
\end{eqnarray*}
for $dg_{\Z}\Vect_k$   the category  of  $\Z$-graded chain complexes over $k$, and $\EExt$ the hyperext functor of \cite{W} \S 10.7.
\end{proposition}
\begin{proof}
Observe that  $\H^*(E_n)$ is cohomology of the complex 
$$
\underline{\Hom}_n(U_{\bt}, \top U_{\bt}) \to \underline{\Hom}_n(U_{\bt}, \top^2 U_{\bt}) \to \ldots
$$
associated to the monad $\top$ (as in \cite{W} \S 8.7),  so for $n>0$, $N^s_n\H^*( E)=\H^*(N^s_nE)$ is cohomology of the complex
$$
\Hom(U_*, (\top U)_*[n]) \to \Hom(U_*, (\top^2 U)_*[n]) \to \ldots.
$$

Now, $(\top U)_*= \bot(U_*)$, and the  augmented cosimplicial complex $\xymatrix@1{U_* \ar[r] & \bot (U_*) \ar@<.5ex>[r] \ar@<-.5ex>[r] &\bot^2 (U_*) \ar[r]& \ldots}$ is canonically contractible (in the sense of \cite{W} 8.4.6), giving
$$
N^s_n\H^i( E)= \left\{ \begin{matrix} \Hom(U_*, U_*[n]) & i=0\\ 0 & i>0, \end{matrix} \right.
$$ 
for $n>0$, so $E$ is quasi-smooth, and $\H^{-n}(\lrcorner E)=\EExt^{-n}_{dg_{\Z}\Vect_k}(U_{\bt},U_{\bt})$ for $n \ge 0$.

Since  $\ker( \gamma_U:\bot U_* \to U_*)= U_*[-1]$,
the cosimplicial normalisation $N_c^n( \bot^{\bt} U_*)= N_c^{n-1} ( \bot^{\bt} U_*)[-1]$, so  $N_c^n( \bot^{\bt} U_*)=U_*[-n]$. Thus   $N_c\CC^{\bt}(E_0)$ is just 
$$
\Hom(U_*, U_*) \xra{\ad_d} \Hom(U_*, U_*[-1]))\xra{\ad_d} \Hom(U_*, U_*[-2]))\xra{\ad_d} \ldots,
$$
so $\H^{i}(\lrcorner E)=\EExt^{i}_{dg_{\Z}\Vect_k}(U_{\bt},U_{\bt})$  for $i> 0$.
\end{proof}

\begin{remarks}\label{moreinterest}
\begin{enumerate}
\item Dually, we may consider deformations of (non-negatively graded) cochain complexes. This simplicial category is monadic over graded modules.

\item We may incorporate the constructions of this section into more interesting examples. For instance, deformations of a complex of $\O_X$-modules on an algebraic space $X$ are given by considering the diagram
$$
\xymatrix@=8ex{
dg\O_X\Mod(X) \ar@<1ex>[r]_{\top} \ar@<-1ex>[d]_{u^{-1}} 
&\ar@<1ex>[l]^{\O_X\ten -} dg\Mod(X)  \ar@<-1ex>[d]_{u^{-1}} 
\\
\ar@<-1ex>[u]_{u_*G}^{\dashv}	g(u^{-1}\O_X)\Mod(X') \ar@<1ex>[r]_{\top} 
&\ar@<1ex>[l]^{ (u^{-1}\O_X)\ten -} \ar@<-1ex>[u]_{u_*G}^{\dashv} g\Mod(X'),  
}
$$
of simplicial categories, where $u:X'\to X$ is the map to $X$ from its set of geometric points.
The resulting pre-DDC will be quasi-smooth whenever $\Ext^i_{\O_X}(M_m,M_n)=0$ for all $i>0$ and $n>m$.
\end{enumerate}
\end{remarks}

\subsection{Simplicial complexes}\label{sections}

\begin{definition} Define $s\FMod(A)$ to be  the category  of simplicial flat modules over $A$. We make this into a simplicial category by setting
$$
\underline{\Hom}(U_{\bt},V_{\bt})_n:= \Hom(\Delta^n\ten U_{\bt},V_{\bt}),
$$
where, for a set $X$ and module $U$, we set $U\ten X:= \bigoplus_{x \in X} U$. 
\end{definition}

\begin{definition}\label{delta*}
Define $\Delta_*$ to be the subcategory of the ordinal number category $\Delta$ containing only  those morphisms fixing $0$.  Given a category $\C$, define the category $s_+\C$ of almost simplicial complexes  (resp.the category $c_+\C$ of  almost cosimplicial complexes) in $\C$ to consist of  functors $(\Delta_*)^{\op} \to \C$ (resp. $\Delta_* \to \C$). Thus an almost simplicial object $X_*$ consists of objects $X_n \in \C$, with all of the operations $\pd_i, \sigma_i$ of a simplicial complex except $\pd_0$, satisfying the usual relations. Similarly, an almost cosimplicial complex has all of the coface and coboundary  operations except $\pd^0$. 
\end{definition}

From now on, we will denote simplicial sets by $X_{\bt}$, and their underlying almost simplicial complexes  by $X_*$.

\begin{definition} Define 
a simplicial structure on the category $s_+\FMod(A)$  by setting
$$
\underline{\Hom}(U_*,V_*)_n:= \Hom(\Delta^n_*\ten U_*,V_*).
$$
\end{definition}

\begin{remark}\label{doldkangraded}
Recall that the Dold-Kan correspondence gives an equivalence $N: s\Mod \to dg\Mod$ of categories,
by the formula $N(V)_n = \bigcap_{i=1}^n \ker(\pd_i:V_n \to V_{n-1})$, with $d:=\pd_0$. Observe that this   extends to an equivalence $N:s_+\Mod \to g\Mod$ of categories, given by the same formula. This is only a weak equivalence of simplicial categories, not an equivalence.
\end{remark}

\begin{lemma}\label{Gdef}
The forgetful functor $s\FMod \to s_+\FMod$ of simplicial categories has a right adjoint $G_{\pd}$, 
and the resulting adjunction is comonadic.
\end{lemma}
\begin{proof} Let
$
(G_{\pd} V_*)_n := V_n\oplus V_{n-1} \oplus \ldots \oplus V_0,
$
with operations
\begin{eqnarray*}
\pd_i(v_n,\ldots,v_0)= (\pd_iv_n,\pd_{i-1}v_{n-1}, \ldots, \pd_1v_{n-i+1}, v_{n-i-1}, \ldots, v_1,v_0)\\
\sigma_i(v_n,\ldots,v_0)= (\sigma_iv_n,\sigma_{i-1}v_{n-1}, \ldots, \sigma_0v_{n-i}, v_{n-i}, \ldots, v_1,v_0).
\end{eqnarray*}
The unit $\alpha:U_{\bt} \to (G_{\pd}U_*)_{\bt}$ of the adjunction is 
$$
\alpha(u)=(u,\pd_0u, \pd_0^{\phantom{0}2}u, \ldots, \pd_0^{\phantom{0}n}u),
$$
 for any simplicial complex $U_{\bt}$, and $u \in U_n$.
 The co-unit $\gamma:G_{\pd}V_* \to V_*$ is the map $\gamma(v_n,\ldots,v_0)=v_n$ of almost simplicial complexes.  
\end{proof}

\begin{remark}\label{suspdef}
The forgetful functor $V_{\pd}$ also has a left adjoint, which does not respect the simplicial structure of the categories. It is given by $\cL_{\pd}(V_*)_n=V_{n+1}$, with $\pd_i^{\cL_{\pd} V}=\pd_{i+1}^V$,  $\sigma_i^{\cL_{\pd} V}=\sigma_{i+1}^V$, and unit $\sigma_0:  V_*\to V_{\pd}\cL_{\pd}(V_*)=\cL_{\pd}(V)_*$.  Note that $\cL_{\pd}V_{\pd}$   is the functor $\mathrm{DEC}^{\op}$ defined on simplicial sets in \cite{glenn}. 
\end{remark}

\begin{definition}\label{xidef}
Define objects $\Xi^n \in s_+\Set$ by $\Xi^n_m=\Hom_{\Delta_*}([m],[n])$, and let $\pd\Xi^n$ be the boundary of $\Xi^n$ (i.e. the union of the images of all maps $\Xi^{n-1} \to \Xi^n$).
Note that $\cL(\Xi^n)=\Delta^n$ and,  for $n>0$, $\cL(\pd\Xi^n)=\L_0^n$, the $0$th horn.
\end{definition}

\begin{lemma}\label{lx0} For $X \in s_+\Set$, $X_0 \to \cL(X)$ is a weak equivalence.
\end{lemma}
\begin{proof}
Any injective map $f:Z \to X$ in $s_+\Set$ is an inductive limit of pushouts of maps $\pd\Xi^n \to \Xi^n$ for $n\ge 0$. If $f_0:Z_0 \to X_0$ is an isomorphism,  we may take $n>0$ only.   Then $\cL(f):\cL(Z) \to\cL(X)$ is an inductive limit of pushouts of maps $\L_0^n \to \Delta^n$, so is a trivial cofibration. Taking $Z=X_0$ gives the required result.
\end{proof}

Let $\bot=V_{\pd}G_{\pd}$ and $\top=G_{\pd}V_{\pd}$, for $V_{\pd}:s\FMod(A)\to s_+\FMod(A)$ the forgetful functor.

\begin{proposition}\label{sdeform}
$s_*\FMod$ has uniformly trivial deformation theory, so given $U_{\bt} \in s\FMod(k)$, we may lift $U_* \in s_+\FMod(k)$
 to $\tilde{U}_* \in s_+\FMod(\L)$. The conditions of Theorem \ref{ddcmain} are satisfied, and the  pre-DDC 
$$
E^n(A):= \underline{\Hom}_{s_+\FMod(A)}(  \tilde{U}_*, \bot^n \tilde{U}_*)_{V(\alpha^n_{U})}
$$
 is then a DDC, with cohomology $\H^*(\lrcorner E)$ given by the complex
$$
 \ldots \xra{\ad_d} \Hom_{g_{\Z}\Mod_k}( NU_*, NU_*[-n]) \xra{\ad_d} \Hom_{g_{\Z}\Mod_k}( NU_*, NU_*[-n-1])\xra{\ad_d} \ldots, 
$$
i.e. $\EExt^*_{dg_{\Z}\Vect_k}(NU_{\bt},NU_{\bt})$.
\end{proposition}
\begin{proof}
This is similar to Proposition \ref{dgdeform}. The only difficulty lies in establishing Definition \ref{q12}.(Q2'):
$$
N_n\CC^i(E)= N_n\underline{\Hom}_{s\FMod(k)}(  U_{\bt}, \top^{i+1} U_{\bt})= \Hom_{s\FMod(k)}(  U_{\bt}\ten (\Delta^n/\L^n_0), \top^{i+1} U_{\bt}).
$$

Now, since $k\ten \L^n_0 \to k\ten \Delta^n$ is a weak equivalence admitting a retraction, $P:=U_{\bt}\ten (\Delta^n/\L^n_0)$ is trivially cofibrant, so is a projective object in the category $s\FMod(k)$. The cosimplicial complex $U^{\bt}$ in $s\FMod(k)$ given by $U^i:=\top^{i+1} U_{\bt}$ is a resolution of $U:=U_{\bt}$ (since the augmented cosimplicial complex $V_{\pd}U \to V_{\pd}U^{\bt}$ is contractible). Thus
$$
\H^*N_n\CC^{\bt}(E)=\EExt^*_{s\FMod(k)}(P,U^{\bt})= \Ext^*_{s\FMod(k)}(P,U)=\Hom_{s\FMod(k)}(P,U).
$$

We have therefore shown that for $n>0$, 
$$
N_n\H^i(E)= \left\{ \begin{matrix} N_n\underline{\Hom}_{s\FMod(k)}(U_{\bt},U_{\bt}) & i=0\\
0 & i>0. \end{matrix} \right.
$$
\end{proof}

\subsection{Simplicial algebras}\label{sectionsalg}

Although the results in this section are expressed for commutative algebras, they will hold for any category equipped with a suitable forgetful functor to flat modules, and in particular algebras over any operad.

\begin{definition}
Let $\FAlg(A)$ be the category of flat (commutative) $A$-algebras, with $s\FAlg(A):= \FAlg(A)^{\Delta^{\op}}$ and $s_+\FAlg(A):= \FAlg(A)^{\Delta_*^{\op}}$ . Recall that, for $K \in \bS$ and $R \in s\FAlg(A)$, we define $R \ten K \in s\FAlg(A)$ by
$$
(R\ten K)_n := R_n^{\ten K_n}= \overbrace{R_n \ten R_n \ten \ldots \ten R_n }^{|K_n|}.
$$
Define $\ten K: s_+\FAlg \to s_+\FAlg$ by the same formula.

Now, we make  $s\FAlg(A), s_+\FAlg(A)$ into simplicial categories by setting
\begin{eqnarray*}
\underline{\Hom}_{s\FAlg}(R_{\bt},S_{\bt})_n &:=& \Hom_{s\FAlg}(\Delta^n\ten R_{\bt},S_{\bt}),\\ 
\underline{\Hom}_{s_+\FAlg}(R_*,S_*)_n&:=& \Hom_{s_+\FAlg}(\Delta^n_*\ten R_*,S_*).
\end{eqnarray*}
\end{definition}

We now consider the commutative diagram
$$
\xymatrix@=8ex{
s\FAlg \ar@<1ex>[r]_{\top}^U \ar@<-1ex>[d]_{V_{\pd}}
&\ar@<1ex>[l]^{\Symm}  s\FMod  \ar@<-1ex>[d]_{V_{\pd}} 
\\
\ar@<-1ex>[u]_{G_{\pd}}^{\dashv}s_+\FAlg	 \ar@<1ex>[r]_{\top}^U 
&\ar@<1ex>[l]^{\Symm} \ar@<-1ex>[u]_{G_{\pd}}^{\dashv} s_+\FMod  
}
$$
of adjunctions of homogeneous simplicial category-valued functors on $\C_{\L}$.

\begin{definition}\label{smodrdef}
Given an almost simplicial $A$-algebra $R_*$, define the category $s_+\Mod(R_*)$ to consist of almost simplicial $A$-modules  $M_*$ equipped with an associative multiplication $R_*\ten M_* \to M_*$, respecting the almost simplicial structures. This has a simplicial model structure (by applying \cite{Hirschhorn} Theorem 11.3.2 to the forgetful functor $s_+\Mod(R_*)\to s_+\Mod(A)$).  

All objects of $s_+\Mod(R_*)$ are fibrant.  Since $A\ten (\pd \Delta^n)_* \to A\ten ( \Delta^n)_*$ has a retraction in $s_+\Mod(A)$, we see that for all cofibrant $C_* \in s_+\Mod(R_*)$,
$$
\underline{\Hom}_{s_+\Mod(R_*)}(C_*,M_*)
$$
is trivially fibrant.
\end{definition}

\begin{proposition}\label{salgdeform}
Fix $R_{\bt} \in s\FAlg(k)$, set $M_*:= UV_{\pd}R_{\bt}$, and choose a lift $\tilde{M}_* \in s_+\FMod(\L)$.

Theorem \ref{ddcmain} gives a pre-DDC
$$
E^n:= \underline{\Hom}_{s_+\FMod}((U\Symm)^n\tilde{M}_*,(V_{\pd}G_{\pd})^n\tilde{R}_*)_{UV_{\pd}(\alpha^n_{R}\circ \vareps^n_R)}.
$$

Moreover,
if  $R_{\bt} \in s\FAlg(k)$     is cofibrant,   then the pre-DDC $E$
of Theorem \ref{ddcmain} is a DDC.  The almost simplicial $k$-algebra $R_*$ then lifts to $ \tilde{R}_* \in s_+\FAlg(\L)$, and $E$ is quasi-isomorphic to the DDC defined by 
$$
(E')_n=\underline{\Hom}_{s_+\FAlg}(\tilde{R}_*,(V_{\pd}G_{\pd})^n\tilde{R}_*)_{V_{\pd}(\alpha^n_{R})}
$$ coming from the comonadic adjunction 
$\xymatrix@1{s\FAlg \ar@<1ex>[r]_{\bot}^{V_{\pd}} & s_+\FAlg \ar@<1ex>[l]^{G_{\pd}}}$.
\end{proposition}
\begin{proof}
It is straightforward to verify   \cite{paper2} equations \ref{paper2-eqnone}--\ref{paper2-eqnfive}, since all our constructions commute with forgetful functors, so $E$ is a pre-DDC.  Since $s_+\FMod$ is uniformly of trivial deformation theory, Definition \ref{q12}.(Q1) is satisfied by $E$. 

 To establish quasi-smoothness, we must compute cohomology groups.
Given a $k$-algebra $S$, recall that the cotangent complex is given by $\bL_n(S/k)=J_n/(J_n)^2$, where $J_n$ is the kernel of the diagonal map $(\Symm U)^{n+1}(S)\ten_kS \to S$.
The cosimplicial complex $\CC^{\bt}(E_n)$ is then given by  $\CC^{m}(E_n)= \Hom_{s_+\Mod(R_*)}( \bL_m(R/k)_*\ten \Delta^n, G_{\pd}^{m}R_*)$. Thus $\H^*(E_n)$ is the total cohomology of the double complex 
$$
\CC^{ij}= \Hom_{s_+\Mod(R_*)}( \bL_i(R/k)_*\ten \Delta^n, G_{\pd}^{j}R_*)= \Hom_{s_+\Mod(R_*)}( \bL_i(R/k)_*, (G_{\pd}^{j}R_*)^{\Delta^n}).
$$ 

Now, if $R_{\bt}$ is cofibrant, the  augmented complex $\bL_{\bt}(R_*) \to \Omega(R_*/k)$ is a levelwise cofibrant resolution in $s_+\Mod(R_*)$. Since all maps in $s_+\Mod(R_*)$ are weak equivalences, cofibrant modules are projective, so the complex is
contractible.

Define Andr\'e-Quillen cohomology on $s_+\Mod(R_*)$ by $D^q(R_*/k,M_*):=\H^q\Hom_{s_+\Mod(R_*)}(\bL_{\bt}(R/k)_*, M_*)$. Given a small extension $A \to B$ with kernel $I$, and a flat almost simplicial $B$-algebra $S_*$,   note that the obstruction to lifting $S_*$ to a flat $A$-algebra lies in $D^2(S_*/B,S_*\ten_BI)=D^2((S_*\ten_Bk)/k, S_*\ten_BI)$, applying \cite{paper2} Theorem \ref{paper2-main} to the adjunction
$$
\xymatrix@1{ s_+\FAlg \ar@<1ex>[r]_{\top} &  s_+\FMod \ar@<1ex>[l]^{\Symm}}.
$$
This ensures that $R_*$ lifts to some $ \tilde{R}_* \in s_+\FAlg(\L)$, so $E'$ can be defined.

Similarly to \cite{paper2} \S \ref{paper2-smooth}, we see that $E'$ is a levelwise smooth DDC, and that $\H^*(E'_n)$ is cohomology of the complex $(\CC')^{m}= \Hom_{s_+\Mod(R_*)}( \Omega(R/k)_*\ten \Delta^n, G_{\pd}^{m}R_*)$. 

Now,  the canonical map $E' \to E$ gives quasi-isomorphisms  $E'_n \to E_n$ for all $n$. We know that $E$ automatically satisfies (Q1). Since $\Omega(R_*/k)$ is a  cofibrant $R_*$-module, the tangent space $\CC(E')^n= \underline{\Hom}_{s_+\Mod(R_*)}(\Omega(R_*/k), G_{\pd}^nR_*)$ is trivially fibrant, so $E'$ also satisfies (Q1).

It only remains to show that $E'$ satisfies (Q2'); the proof of Proposition \ref{sdeform} adapts. 
\end{proof}

\begin{remark}
We may weaken the condition that $R_{\bt}$ be cofibrant to requiring that the cotangent complex $\diag\bL_{\bt}(R_{\bt}/k)$ of $R_{\bt}$ is equivalent in $\Ho(s\Mod(R_{\bt}))$ to $\Omega(R_{\bt}/k)$, and that the latter is cofibrant. If a $k$-algebra $R$ (with constant simplicial structure) is  smooth, \cite{Ill1} Proposition III.3.1.2 implies that this holds.
\end{remark}

\begin{definition}
Given a simplicial $k$-algebra $R_{\bt}$, and a simplicial $R_{\bt}$-module $M_{\bt}$, define the simplicial vector space $\underline{\Der}(R,M)$ of derived derivations by 
$$
\underline{\Der}(R_{\bt},M_{\bt})_n := \Der(R_{\bt}\ten \Delta^n, M_{\bt}),
$$
the set of simplicial $k$-algebra morphisms $f: R_{\bt}\ten \Delta^n \to R_{\bt}\oplus M_{\bt}\eps$ extending the canonical map $ R_{\bt}\ten \Delta^n \to R_{\bt}$, where $\eps^2=0$.
\end{definition}

\begin{remark}\label{salgcoho}
For $R_{\bt}\in s\Alg(k)$ cofibrant, and $E$ as in Proposition \ref{salgdeform}, $\H^n(E)=\pi_{-n}\underline{\Der}_{s\Alg}(R_{\bt},R_{\bt})$ for $n \le 0$. For $n>0$, $\H^n(E)$ is the $n$th cohomology of the cosimplicial complex
$$
C^n:= \Der_{s\Alg}(\cL_{\pd}^{n+1} R_{\bt},R_{\bt})
$$
associated to the comonad $\cL_{\pd}$ of Definition \ref{suspdef}. 
\end{remark}

\begin{proposition}\label{sconsistent}
If $R_{\bt} \to R$ is a cofibrant resolution of a $k$-algebra $R$, then  the DDC $E$ of Proposition \ref{salgdeform} is quasi-isomorphic to $DF$, for $F$ the SDC
$$
F^n=\Hom_{\Mod}(\Symm^n\tilde{M},\tilde{M})
$$ 
from \cite{paper1} \S\ref{paper1-alg}, and $\tilde{M} \in \FMod(\L)$ lifting the $k$-module $M$ underlying $R$. 
\end{proposition}
\begin{proof}
By Proposition \ref{weakin}, we may assume that $R_{\bt}$ is the standard resolution $R_n=\bot^{n+1}R$, with $\tilde{R}_n= \Symm^{n+1}(\tilde{M})$.

Then we have $E$ quasi-isomorphic to the DDC $E'$ given by 
$$
(E')^n_K= \Hom_{s_+\Alg}(K_*\ten \tilde{R}_*, G_{\pd}^n\tilde{R}_*)=\Hom_{s_+\Mod}(K_*\ten\Symm^*\tilde{R},G_{\pd}^n\tilde{R}_*).
$$
 
The augmentation $\vareps:\tilde{R}_*\to \tilde{R}$ in $s_+\Mod(\L)$ gives us a map
\begin{eqnarray*}
\chi:(E')^n_K\to \Hom_{s_+\Mod}(K_*\ten\Symm^*\tilde{R},G_{\pd}^n\tilde{R})&=& \Hom_{s_+\Mod}(\cL^n(K_*\ten\Symm^*\tilde{R}),\tilde{R})\\
&=&\Hom_{\Mod}(K_n\ten \Symm^n\tilde{R},\tilde{R}).
\end{eqnarray*}
But this is just $(DF)_K^n$, and it is straightforward to check that $\chi:E'\to DF$  respects all the SDC operations. 

Since $R_{\bt} \to R$ is a resolution, we get a weak equivalence
$$
\z^n\CC^{*}(E)=\z^n\underline{\Der}_{s\Alg}(R_{\bt},G^{*+1}R_{\bt})\to \z^n\underline{\Der}_{s\Alg}(R_{\bt},G^{*+1}R)
$$
Now by Lemma \ref{lx0},
\begin{eqnarray*}
\underline{\Der}_{s\Alg}(R_{\bt},G^{n+1}R)_m&=& \Der_{\Alg}( \pi_0\cL^{n+1}(\Delta^m\ten R_{\bt}),R)\\
&=&\Der_{\Alg}(R_n,R)^{\Delta^m_n}\\
&=&\Der_{\Alg}(\Symm^nR,R)^{\Delta^m_n}\\
&=& \CC^{n}(DF)_m.
\end{eqnarray*}
Since $\CC^{\bt}(E)$ is quasi-smooth,  $\H^{n-i}(\lrcorner E)=\pi_i\z^n\CC^{*}(E)$, so 
$$
\H^*(\lrcorner E)\cong \H^*(\lrcorner DF )=\H^*(F),
$$
which
 is just Andr\'e-Quillen cohomology $\bD^*_k(R,R)$, and so $\chi$ is a quasi-isomorphism.
\end{proof}

\begin{remark}
Propositions \ref{salgdeform} and \ref{sconsistent} together imply that derived deformations of a $k$-algebra $R$ are equivalent to  derived deformations of the operation $\pd_0$ on any cofibrant resolution $R_{\bt} \to R$.
\end{remark}

\section{Deformations of Artin stacks and simplicial schemes}\label{stacksn}

The problem we now wish  consider is that of deforming of an algebraic  stack $\fX$. We may take a smooth simplicial hypercovering $X_{\bt} \to \fX$, with each $X_n$ a  disjoint union of affine schemes (similarly to the proof of \cite{olssartin} Theorem 11.1), and our first step will be  to consider derived deformations of $X_{\bt}$. 

\subsection{Cosimplicial algebras}\label{calgsn}

Let
$X_{\bt}$ be a simplicial affine scheme.  Equivalently, we may consider the cosimplicial algebra $[n] \mapsto \Gamma(X_n, \O_{X_n})$. 
\begin{definition}\label{csimpstr}
The categories  $c\FMod(A)$, $c\FAlg(A)$,  $c_+\FMod(A)$, and $c_+\FAlg(A)$ (as given in Definition \ref{delta*}) can be made into simplicial categories (i.e. enriched in simplicial sets) by setting
$
(S^K)^n:= (S^n)^{K_n}
$
 for $K \in \bS$, with structure maps $(S^K)(f)= S(f)^{K_n} \circ K(f)^*: (S^m)^{K_m} \to (S^n)^{K_n}$, for  morphisms $f$ in $\Delta$. We then define the simplicial $\HHom$ functor by 
$$
\underline{\Hom}(R,S)_n:= \Hom(R, S^{\Delta^n}).
$$
\end{definition}

There is the following diagram of monadic adjunctions of functors $\C_{\L} \to s\Cat$:
$$
\xymatrix@=8ex{
c\FAlg(A) \ar@<1ex>[r]^-{U_{\alg}}_-{\top} \ar@<-1ex>[d]_{U_{\pd}}^{\vdash} 
&\ar@<1ex>[l]^-{ \Symm }   c\FMod(A)    \ar@<-1ex>[d]_{U_{\pd}}^{\vdash} 
\\
\ar@<-1ex>[u]_{F_{\pd}}c_+\FAlg(A)	 \ar@<1ex>[r]^-{U_{\alg}}_-{\top} 
&\ar@<1ex>[l]^-{\Symm} \ar@<-1ex>[u]_{F_{\pd}}  c_+\FMod(A), 
}
$$
 where $F_{\pd}:c_+\FMod(A)\to c\FMod(A)$ is left adjoint to the forgetful functor $U_{\pd}$, given by 
$$
(F_{\pd} V^*)^n = V^n\oplus V^{n-1} \oplus \ldots \oplus V^0,
$$
with operations dual to those in Lemma \ref{Gdef}. Similarly, $F_{\pd}:c_+\FAlg(A)\to c\FAlg(A)$ is the left adjoint given by 
$$
(F_{\pd} R^*)^n = R^n\ten R^{n-1} \ten \ldots \ten R^0.
$$
The diagram satisfies the following commutativity conditions:
$$
U_{\pd}U_{\alg}=U_{\alg}U_{\pd} \quad \Symm F_{\pd}=F_{\pd}\Symm, \quad U_{\pd}\Symm=\Symm U_{\pd}.
$$
These adjunctions combine to give a monadic adjunction
$$
\xymatrix@1{ c\FAlg(A) \ar@<1ex>[r]^-{U_{\pd} U_{\alg}}_-{\top}  
&\ar@<1ex>[l]^-{\Symm  F_{\pd}} c_+\FMod(A)}.
$$

\begin{lemma}
$c_+\FMod(A)$ has uniformly trivial deformation theory.
\end{lemma}
\begin{proof}
This is essentially the same as Proposition \ref{sdeform}.
\end{proof}

\begin{proposition}\label{stackprop1}
In the scenario above,  Theorem \ref{ddcmain} gives a pre-DDC
$$
E^n:= \underline{\Hom}_{c_+\FMod}(\top^n\widetilde{U_{\pd}U_{\alg} R},\widetilde{U_{\pd}U_{\alg}R}),
$$
 satisfying  Definition \ref{q12}.(Q1), where $\top=U_{\pd}U_{\alg}\Symm F_{\pd}$.
\end{proposition}

We now seek conditions under which the pre-DDC $E$ (or similarly a pre-DDC $E(\bD)$ associated to a diagram as in \S \ref{weakinvart}) is quasi-smooth.

\begin{definition}\label{cmodrdef}
Given a cosimplicial (resp. almost cosimplicial) $A$-algebra $R$, define the category $c\Mod(R)$ (resp.  $c_+\Mod(R)$) to consist of cosimplicial (resp. almost cosimplicial) $A$-modules $M$  equipped with an associative multiplication $R\ten M \to M$, respecting the cosimplicial (resp. almost cosimplicial) structures. These categories have  simplicial structures, with 
$
(M^K)^n:= (M^n)^{K_n},
$   
for $K \in \bS$, the $R$-module structure on $M^K$ coming from the map $R \to R^K$. As usual, denote the left adjoint to $M \mapsto M^K$ by $N \mapsto N\ten K$.
\end{definition}

Given $M \in c\Mod(R)$ and an injective map $K\into L$ in $\bS$, set $M\ten (L/K):= \coker(M\ten K \to M\ten L)$ and $M^{L/K}:= \ker(M^L \to M^K)$.

\begin{definition}
Let $\bot_{\alg}= \Symm U_{\alg}$, $\bot_{\pd}=F_{\pd}U_{\pd}$ and    $\bot=  \Symm  F_{\pd}  U_{\pd} U_{\alg}=\bot_{\pd}\bot_{\alg}$. Given $R \in c\Alg(k)$, define $\bL^{\bot}_n(R) \in  c\Mod(R)$  by the property that 
$$
\Hom_{c\Mod(R)}( \bL^{\bot}_n(R), M^{\bt})\cong\Der_k( \bot^{n+1}R, M^{\bt})
$$
functorial in $M^{\bt} \in c\Mod(R)$. Here, $\Der_k(S^{\bt}, M^{\bt})$ is the set of  morphisms $f: S^{\bt} \to S^{\bt} \oplus M^{\bt}\eps$  in $c\Alg_k$ extending the identity, where $\eps^2=0$.

Define $\bL_n(R) \in c\Mod(R)$ by 
$$
\Hom_{c\Mod(R)}( \bL_n(R), M^{\bt})\cong\Der_k( \bot_{\alg}^{n+1}R, M^{\bt}).
$$

Observe that $\bL_{\bt}(R)$ and $\bL_{\bt}^{\bot}(R)$ both form simplicial complexes in $c\Mod(R)$.
\end{definition}

\begin{definition}
Given an object $R$ of   $c\FAlg(A)$ (resp. $c_+\FAlg(A)$), we may extend $R$ uniquely to a cocontinuous functor $R: \bS \to \FAlg(A)$ (resp. $R: s_+\Set \to \FAlg(A)$) extending the functor $R: \Delta \to \FAlg(A)$ (resp. $R: \Delta_+ \to \FAlg(A)$)
given by  $R(\Delta^n)=R^n$ (resp. $R(\Xi^n)=R^n$, for $\Xi$ as  in Definition \ref{xidef}). 
\end{definition}

\begin{lemma}\label{lalg}
For all $m$, the simplicial complex $\bL^{\bot}_{\bt}(R)^m$ is a model for the cotangent complex of $R^m$.
\end{lemma}
\begin{proof}
Write $\bot_nR:= \bot^{n+1}R$; these form a simplicial complex $\bot_{\bt}R$ in $c\Alg(k)$.
We need to show that  $(\bot_{\bt}R)^m$ is a cofibrant resolution of $R^m$ in $s\Alg(k)$.
If we apply the forgetful functor $U_{\pd}U_{\alg}$ to the augmented simplicial complex $\bot_{\bt}R\to R $, we see that it becomes contractible. In particular, this implies that $\bot_{\bt}R \to R$ is contractible as an augmented complex of $k$-vector spaces, so it is a resolution. 
\end{proof}

\begin{lemma}\label{lproj}
$\bL^{\bot}_n(R)$ is a projective object of $c\Mod(R)$, and $U_{\pd}\bL_n(R), U_{\pd}\bL^{\bot}_n(R)$ are both projective objects of $c_+\Mod(R^*)$.
\end{lemma}
\begin{proof}
By adjointness,  
$$
\Der_k( \bot^{n+1}R,M^{\bt}) \cong \Hom_{c_+\Mod_k}(U_{\pd}U_{\alg} \bot^nR, M^*),  
$$
so  $\Der_k(\bot^{n+1}R, - )$ defines a right exact functor, hence  $\bL^{\bot}_n(R)$ is  projective. The other results follow similarly.
\end{proof}

\begin{lemma}\label{swapbot}
There is a natural transformation $F_{\pd}U_{\alg} \to U_{\alg}F_{\pd}$, giving  transformations $\bot_{\pd}\bot_{\alg} \to \bot_{\alg}\bot_{\pd}$. 
\end{lemma}
\begin{proof}
The transformation is given on level $n$ by 
$$
R^0 \oplus R^1 \oplus \ldots \oplus R^n \ni \sum_{i=0}^n r_i \mapsto \sum_{i=0}^n 1\ten\ldots \ten 1 \ten r_i \ten 1 \ten \ldots 1 \in  R^0 \ten R^1 \ten \ldots \ten R^n.
$$
\end{proof}

\begin{definition}\label{gsmooth}
A morphism $f:X_{\bt}\to Y_{\bt} $ of simplicial schemes  over $A$ is said to be  quasi-smooth (resp. trivially smooth) if  the morphism
$$
\Hom_{\bS}(L,X_{\bt})\to \Hom_{\bS}(K,X_{\bt})\by_{\Hom_{\bS}(K,Y_{\bt})}\Hom_{\bS}(L,Y_{\bt})
$$
of affine schemes is 
smooth 
for all trivial cofibrations (resp. all cofibrations) $K \to L$ of finite simplicial sets. The map $f$ is said to be smooth if it is quasi-smooth and $f_0:X_0 \to Y_0$ is smooth. 

 We say that a morphism $R^{\bt}\to S^{\bt}$ in  $c\Alg(A)$ is   quasi-smooth (resp.  trivially smooth, resp. smooth) if $\Spec S^{\bt} \to \Spec R^{\bt}$ is so.
\end{definition}

\begin{lemma}\label{gsmoothslack}
In Definition \ref{gsmooth}, we may replace cofibrations (resp. trivial cofibrations) $K \to L$ by generating cofibrations $\pd\Delta^n \to \Delta^n$ (resp.  generating trivial cofibrations $\L^n_k \to \Delta^n$). 
\end{lemma}
\begin{proof}
This follows because every cofibration (resp. trivial cofibration) is a composition of pushouts of generating cofibrations (resp.  generating trivial cofibrations), and the fact that smooth morphisms are closed under pullback and finite composition.
\end{proof}

\begin{lemma}\label{smoothcdns}
A morphism $f:X_{\bt}\to Y_{\bt} $ of simplicial  schemes is quasi-smooth (resp.  trivially smooth, resp. smooth) if and only if the following conditions hold:
\begin{enumerate}
\item for all square-zero extensions $A \onto B$ of $k$-algebras, the map $X_{\bt}(A) \to X_{\bt}(B)\by_{Y_{\bt}(B)}Y_{\bt}(A) $ is a fibration (resp. a trivial fibration, resp. a surjective fibration) in $\bS$.


\item for all  all vertices $v\in \Delta^n_0$ the maps $v^*:X_n \to X_0$ (resp. the schemes $X_n$, resp. the schemes $X_n$) are locally of finite presentation. 
\end{enumerate}
\end{lemma}
\begin{proof}
This follows from the fact that a morphism is smooth if and only if it is quasi-smooth and locally of finite presentation, and that $U \to V$ is locally of finite presentation if and only if
the map $U(A_{\alpha})\to U(\varinjlim  A_{\alpha})\by_{V(\varinjlim  A_{\alpha})}\varinjlim V(A_{\alpha})$ is an isomorphism. We also use the result that if $g \circ f$ is locally of finite presentation, then  $f$ must also be so.
\end{proof}

\begin{corollary}\label{powerqs}
For all cofibrations  $i:K \to L$ of finite simplicial sets,  and $f:X\to Y $ a quasi-smooth morphism of simplicial affine schemes, the map
$$
g:X^L \to X^K\by_{Y^K}Y^L
$$
is quasi-smooth. Moreover, if either  $i$ or $f$ is trivial, then so is $g$.
\end{corollary}

\begin{lemma}\label{omegaproj}
If $R\to S$ is  a trivially smooth map in $c\Alg$, then $\Omega(S/R)$ is projective in $c\Mod(S)$.
\end{lemma}
\begin{proof}
By definition, we know that $R(L)\ten_{R(K)}S(K) \to S(L)$ is smooth for all cofibrations $K \into L$ of finite simplicial sets.

Given   $M^{\bt}  \in c\Mod(S)$ and $K \in \bS$, define $M(K) \in \Mod(S(K))$ by $S(K) \oplus M(K)\eps=(S\oplus M\eps)(K)$. Note that  $\Omega(S/R)(K)=\Omega(S(K)/R(K))$. 

Take a surjection $L^{\bt} \to N^{\bt}$ in $c\Mod(S)$ and a morphism $f:\Omega(S/R) \to N^{\bt}$. We will construct a lifting $\tilde{f}$ of $f$ inductively.  Assume that we have $R^m$-linear maps   $\tilde{f}^m: \Omega(S/R)^m \to L^m$ lifting $f$ compatibly with the cosimplicial operations, for all $m<n$. If $M^nL$ denotes the $m$th matching object (as in \cite{sht} Lemma VII.4.9), then extending $\tilde{f}$ compatibly to $\Omega(S/R)^n$ amounts to finding a lift
$$
\xymatrix{
\Omega(S/R) (\pd \Delta^n)\ten_{S(\pd\Delta^n)}S^n \ar[r] \ar[d]& L^n\ar[d]^{\alpha}\\
\Omega(S/R)^n \ar[r]\ar@{-->}[ur] & N^n\by_{M^{n-1}N}M^{n-1}L.
}
$$

Now, $T:=R^n\ten_{R(\pd\Delta^n)} S(\pd\Delta^n) \to S^n$ is smooth, as is $R(K) \to S(K)$ for all $K$. Thus the sequence
$$
0 \to \Omega(S(\pd\Delta^n)/ R(\pd\Delta^n))\ten_{S(\pd\Delta^n)}S^n \to \Omega(S^n/R^n) \to \Omega(S^n/T) \to 0
$$
is exact, with all terms projective. Since $\alpha$ is surjective, projectivity of $\Omega(S^n/T)$ gives the required lift.
\end{proof}

\begin{lemma}
If $R$ is quasi-smooth, then for all trivial cofibrations $K \to L$ of finite simplicial sets, $\Omega(R)\ten (L/K)$ is projective in $c\Mod(R)$.
\end{lemma}
\begin{proof}
By Corollary \ref{powerqs}, $R\ten K \to R\ten L$ is trivially smooth. Since $\Omega(R)\ten X=\Omega(R\ten X)\ten_{R\ten X}R$,  projectivity follows from Lemma \ref{omegaproj}.  
\end{proof}

\begin{definition}
Let $\cL^{\vee}: c_+\Mod\to c\Mod$ be right adjoint to $U_{\pd}$.
\end{definition}

\begin{lemma}\label{losebot} 
For $R \in c\Alg(k)$  and an injective map $f:Z \to X$ in $s_+\Set$, 
there is an isomorphism
$$
 \H^i\Hom_{c_+\Mod(R^*)}(\bL_{\bt}(R^*), N \ten k^{X/ Z}) \cong \EExt^i_{c\Mod(R)}(\bL^{\bot}_{\bt}(R)\ten(\cL X/\cL Z),\cL^{\vee} N)
$$
for all $N \in c_+\Mod(R^*)$.
\end{lemma}
\begin{proof}
By Lemma \ref{lproj}, $\bL^{\bot}_{n}(R)$ is projective, so  $\bL^{\bot}_{n}(R)\ten(\cL X/\cL Z)$ must also be projective, as  $M \mapsto M^K$ sends surjections to surjections.
  Thus 
\begin{eqnarray*}
\EExt^*_{c\Mod(R)}(\bL^{\bot}_{\bt}(R)\ten(\cL X/\cL Z),\cL^{\vee} N)&=& \H^*\Hom_{c\Mod(R)}(\bL^{\bot}_{\bt}(R)\ten(\cL X/\cL Z),\cL^{\vee} N)\\
&=& \H^*\Hom_{c\Mod(R)}(\bL^{\bot}_{\bt}(R),\cL^{\vee}N \ten \cL^{\vee}k^{X/ Z})\\
&=& \H^*\Hom_{c_+\Mod(R^*)}( U_{\pd}\bL^{\bot}_{\bt}(R), N \ten k^{X/ Z})\\
&=& \EExt^*_{c_+\Mod(R^*)}( U_{\pd}\bL^{\bot}_{\bt}(R), N \ten k^{X/ Z}).
\end{eqnarray*}

Now, Lemma \ref{swapbot} gives compatible transformations $\bot_{\pd}^{n+1}\bot_{\alg}^{n+1}(R) \to \bot^{N+1}R$. The unit of the adjunction $F_{\pd} \dashv U_{\pd}$ gives compatible transformations  $U_{\pd} \to U_{\pd}\bot_{\pd}^{n+1}$, so there is a map $\bL_{\bt}(R^*)=U_{\pd}\bL_{\bt}(R^{\bt}) \to  U_{\pd}\bL^{\bot}_{\bt}(R)$, which is an equivalence in the derived category by Lemma \ref{lalg}. Hence
\begin{eqnarray*}
\EExt^*_{c_+\Mod(R^*)}( U_{\pd}\bL^{\bot}_{\bt}(R), N \ten k^{X/ Z}) &\cong&\EExt^*_{c_+\Mod(R^*)}( \bL_{\bt}(R^*), N \ten k^{X/ Z})\\
&=& \H^*\Hom_{c_+\Mod(R^*)}(\bL_{\bt}(R^*), N \ten k^{X/ Z}).
\end{eqnarray*}
\end{proof}

\begin{lemma}\label{qskey}
For $R \in c\Alg(k)$  quasi-smooth, there is an exact sequence
$$
0 \to \bL_{\bt}(R^0)\ten_{R^0}R^*\to \bL_{\bt}(R^*)\to \Omega(R^*/R^0)\to 0 
$$
in the derived category of  projective complexes in $c_+\Mod(R^*)$, where the morphism $R^0 \to R^*$ is given in level $n$ by $(\pd^1)^n$.
\end{lemma}
\begin{proof}
There is an  exact sequence
$$
0 \to \bL_{\bt}(R^0)\ten_{R^0}R^*\to \bL_{\bt}(R^*)\to \bL(R^*/R^0)\to 0 
$$
in the derived category.
Since $R^{\bt}$ is quasi-smooth,  the maps $(\pd^1)^n:R^0 \to R^n$ are all smooth, giving $\bL(R^*/R^0)\sim \Omega(R^*/R^0)$
\end{proof}

\begin{proposition}\label{stackqsgood}
 If $R$ is a quasi-smooth object of $c\Alg(k)$, then every   morphism $\rho:R \to S$ in $c\Alg$  is quasi-smooth over $c_+\Mod$, in the sense of Definition \ref{qsmorphism}.
The $\Ext$-groups are then given by
$$
\Ext^i_{c\Alg/c_+\Mod}(\rho):= \left\{ \begin{matrix} \EExt^i_{c\Mod(R)}( \bL_{\bt}^{R/k}, S) & i>0 \\ \pi_{-i}\Hom_{c\Mod(R)}(\Omega_{(R\ten \Delta^{\bt})/k}\ten_{R\ten \Delta^{\bt}}R,S) & i \le 0,  \end{matrix} \right.
$$
where $\Hom_{c\Mod(R)}(\Omega_{(R\ten \Delta^{\bt})/k}\ten_{R\ten \Delta^{\bt}}R,S)$ is the simplicial complex given in level $n$ by $\Hom_{c\Mod(R)}(\Omega_{(R\ten \Delta^{n})/k}\ten_{R\ten \Delta^{n}}R,S)$.
\end{proposition}
\begin{proof}
First observe that since $c_+\Mod$ has uniformly trivial deformation theory, $\rho$ is quasi-smooth over $c_+\Mod$ whenever it is Q2 over $c_+\Mod$. Now,  
$$
\CC^{\bt}_{c\Alg/c_+\Mod} (\rho)= \HHom_{c\Mod(R)}(\bL^{\bot}_{\bt}(R),S),
$$
so for $K \in \bS$,
\begin{eqnarray*}
\underline{\Ext}^*_{c\Alg/c_+\Mod}(\rho)&=& \H^*\HHom_{c\Mod(R)}(\bL^{\bot}_{\bt}(R),S),\\
\underline{\Ext}^*_{c\Alg/c_+\Mod}(\rho)_K &=& \H^*\Hom_{c\Mod(R)}(\bL^{\bot}_{\bt}(R)\ten K,S),\\
&\cong& \EExt^*_{c\Mod(R)}(\bL^{\bot}_{\bt}(R)\ten K,S),
\end{eqnarray*}
the latter isomorphism following since $\bL^{\bot}_{n}(R)\ten K $  is projective.

Now, consider the monad $\top_{\pd}:= \cL^{\vee}U_{\pd}$ on $c\Mod(R)$, and observe that the augmented cosimplicial complex $\top_{\pd}^{\bt+1}M$  given in level $n$ by $\top_{\pd}^{n+1}M$ is a resolution in $c\Mod(R)$, since it becomes contractible on applying $U_{\pd}$. Thus
$$
\underline{\Ext}^*_{c\Alg/c_+\Mod}(\rho)_K \cong \EExt^*_{c\Mod(R)}(\bL^{\bot}_{\bt}(R)\ten K,\top_{\pd}^{\bt+1}S).
$$

Given $f:Z \to X$ in $s_+\Set$ with $f_0$ an isomorphism, by Lemmas \ref{losebot} and \ref{qskey}, we have
$$
\EExt^i_{c\Mod(R)}(\bL^{\bot}_{\bt}(R)\ten (\cL X/\cL Z),\cL^{\vee} N)= \H^i\Hom_{c_+\Mod(R^*)} ( \bL_{\bt}(R^0)\ten_{R^0}R^*,N \ten k^{X/ Z})
$$
for all $i>0$. However,
$$
\Hom_{c_+\Mod(R^*)} ( \bL_{\bt}(R^0)\ten_{R^0}R^*,N \ten k^{X/ Z})=\Hom_{\Mod(R^0)} ( \bL_{\bt}(R^0), (N \ten k^{X/ Z})^0)=0, 
$$
since $(k^{X/ Z})^0=0$.

Hence $\EExt^i_{c\Mod(R)}(\bL^{\bot}_{\bt}(R)\ten (\cL X/\cL Z),\cL^{\vee}N)=0$ for all $i>0$, so the spectral sequence associated to $\top_{\pd}^{\bt+1}$ gives
$$
\EExt^*_{c\Mod(R)}(\bL^{\bot}_{\bt}(R)\ten (\cL X/\cL Z),M)\cong \H^*\Hom_{c\Mod(R)}(\Omega(R)\ten (\cL X/\cL Z), \top_{\pd}^{\bt+1}M ).
$$

Since $\Omega(R)\ten (\cL X/\cL Z)$ is projective (by Lemma \ref{omegaproj}), this is just
\begin{eqnarray*}
\EExt^*_{c\Mod(R)}(\Omega(R)\ten (\cL X/\cL Z), \top_{\pd}^{\bt+1}M )&\cong& \Ext^*_{c\Mod(R)}(\Omega(R)\ten (\cL X/\cL Z), M )\\
&=& \Hom_{c\Mod(R)}(\Omega(R\ten \cL X/R\ten\cL Z)\ten_{R\ten \cL X}R , M).
\end{eqnarray*}

Taking $Z=\pd\Xi^n, X=\Xi^n$, we have $\cL Z=\L^n_0, \cL X=\Delta^n$, and 
$$
N_n\underline{\Ext}^*_{c\Alg/c_+\Mod}(\rho)= \Hom_{c\Mod(R)}(\Omega(R\ten \Delta^n/R\ten \L^n_0)\ten_{R\ten \Delta^n}R , S),
$$
so $\underline{\Ext}^i_{c\Alg/c_+\Mod}(\rho)$ is constant for $i>0$.

Thus $\rho$ is Q2 over $c_+\Mod$, as required. The description of positive $\Ext$-groups follows from Lemma \ref{losebot}, while that of  non-positive $\Ext$-groups follows from the definition of $\pi_{-i}\H^0\CC^{\bt}_{c\Alg/c_+\Mod} (\rho)$.
\end{proof}

\begin{corollary}
For any diagram in $c\Alg$ with quasi-smooth objects, the associated pre-DDC given by Definition \ref{sdcdiagram} and Proposition \ref{enrich} applied to the adjunction 
$$
\xymatrix@1{ c\FAlg(A) \ar@<1ex>[r]^-{U_{\pd} U_{\alg}}_-{\top}  
&\ar@<1ex>[l]^-{\Symm  F_{\pd}} c_+\FMod(A)}
$$
is a DDC by Lemma \ref{q2diagram}, and governs deformations in the simplicial category $c\Alg$ (by Proposition \ref{governdiagram}).
\end{corollary}

\subsubsection{Comparison with deformations of schemes}


In \cite{paper2} \S \ref{paper2-sepn}, an SDC was constructed to describe deformations of a  separated scheme $X$, and  we now wish to compare it with the  DDC above.

Take an open affine cover $(X_{\alpha})_{\alpha \in I}$ of $X$, and set $\check{X}:=\coprod_{\alpha \in I}X_{\alpha}$. Define the simplicial scheme $Z_{\bt}$ by $Z= \cosk_0(\check{X}/X)$, i.e.
$$
Z_n= \overbrace{\check{X}\by_X \check{X}\by_X \ldots \by_X     \check{X}}^{n+1},
$$ 
with $r_n:Z_n \to X$, and $s_n:Z_n\to \check{X}$ given by projection onto the first factor.

The map $v:\check{X} \to X$ gives adjoint functors $v^{-1} \dashv v_*$ on sheaves. 
This yields the following diagram of $\Cat$-valued functors: 
$$
\xymatrix@C=12ex@R=8ex{
\FAlg_A(X) \ar@<1ex>[r]_{\top} \ar@<-1ex>[d]_{v^{-1}}
&\ar@<1ex>[l]^{\Symm_{A}} \FMod_A(X)  \ar@<-1ex>[d]_{v^{-1}}
\\
\ar@<-1ex>[u]_{v_*}^{\dashv}	\FAlg_A(\check{X}) \ar@<1ex>[r]_{\top} 
&\ar@<1ex>[l]^{\Symm_{A}} \ar@<-1ex>[u]_{v_*}^{\dashv} \FMod_A(\check{X}), 
}
$$
where $\FMod_A(Y)$ and  $\FAlg_A(Y)$ denote sheaves of flat $A$-modules and of flat $A$-algebras on $Y$.

\begin{definition}
The SDC $\check{E}^{\bullet}$  of \cite{paper2} \S \ref{paper2-sepn} was then given by 
$$
\check{E}^n(A)=\Hom_{\FMod_{A}(\check{X})}((\Symm_{A})^n \sN\ten A, (v^{-1}v_*)^n \sN\ten A)_{v^{-1}(\alpha^n\circ \vareps^n)},
$$
for $\sN$ a flat $\mu$-adic $\L$-module on $\check{X}$  lifting $v^{-1}\O_X$, with $\alpha^n:\O_X \to (v_*v^{-1})^n\O_X$ coming from the unit of the adjunction, and similarly $\vareps^n:(\Symm_k)^n\O_X \to \O_X$.
\end{definition}

\begin{definition}
Define  functors $\check{\CC}^{\bt}:\FMod_A(X)\to c\FMod(A)$, $\check{\CC}^*:\FMod_A(\check{X})\to c_{+}\FMod(A)$ by 
$$
\check{\CC}^n(\sF):= \Gamma(Z_n, {r_n}^{-1}\sF), \quad \check{\CC}^n(\sG):= \Gamma(Z_n, {s_n}^{-1}\sG),
$$ with the standard cosimplicial operations.
\end{definition}

\begin{lemma}
There are canonical isomorphisms
$$
\CC^*(v^{-1}\sF)=U_{\pd}\CC^{\bt}(\sF) \quad \CC^{\bt}(v_*\sG)\cong \cL^{\vee}\CC^*(\sG).
$$
\end{lemma}

\begin{lemma}
There is a canonical natural transformation $\Symm\circ \check{C}^{\bt} \to \check{C}^{\bt}\circ \Symm$.
\end{lemma}

\begin{proposition}
The SDC $\check{E}$ is quasi-isomorphic to the DDC $E$ of Proposition \ref{stackprop1}, in the sense that $\Def(\check{E})$ and $\Def(E)$ are weakly equivalent (equivalently, $D\check{E}$ and $E$ are quasi-isomorphic DDCs). 
\end{proposition}
\begin{proof}
We have maps
\begin{eqnarray*}
\Hom_{\Mod(\check{X})}(\top_{\alg}^n\sM, (v^{-1}v_*)^n\sM)&\to& \Hom_{c_+\Mod}(\check{\CC}^* (\top_{\alg}^n\sM), \check{\CC}^*(v^{-1}v_*)^n\sM)\\
&=&\Hom_{c_+\Mod}(\check{\CC}^* (\top_{\alg}^n\sM), (\cL^{\vee})^n\check{\CC}^*\sM)\\
&=&  \Hom_{c_+\Mod}(\top_{\pd}^n\check{\CC}^* (\top_{\alg}^n\sM), \check{\CC}^*\sM)\\
&\to& \Hom_{c_+\Mod}(\top_{\alg}^n\top_{\pd}^n\check{\CC}^* (\sM), \check{\CC}^*\sM)\\
&\to&\Hom_{c_+\Mod}(\top^n\check{\CC}^* (\sM), \check{\CC}^*\sM).
\end{eqnarray*}

These are compatible with the SDC operations, giving a morphism $\check{E} \to E_0$ of SDCs. Now, as in Proposition \ref{stackqsgood}, 
$$
\H^*(E_0)= \EExt^*_{\O_{Z_{\bt}}}( \bL^{Z/k}_{\bt}, \O_{Z}).
$$
However, since the maps $r_n:Z_n \to X$ are all open, and hence \'etale, $\bL^{Z/k}_{\bt}$ is quasi-isomorphic to $r^*\bL^X_{\bt}$. Thus
$$
\EExt^*_{Z_{\bt}}( \bL^{Z/k}_{\bt}, \O_{Z})= \EExt^*_{\O_{X}}( \bL^{X/k}_{\bt}, r_*\O_{Z_{\bt}})= \EExt^*_{\O_{X}}( \bL^{X/k}_{\bt}, \O_X),
$$
since $r_*\O_{Z_{\bt}}=r_*r^{-1}\O_X$ is a resolution of $\O_X$. This means that $\check{E} \to E_0$ is a quasi-isomorphism of SDCs.

Finally, to see that $\Def(E_0) \to \Def(E)$ is a quasi-isomorphism, apply Lemma \ref{levelwiseqscoho},  noting that the strictly positive cohomology groups automatically agree. For $n\le 0$,
$$
\H^nE= \H_{-n}\Hom_{\O_Z}(i^*\Omega_{Z^{\Delta^{\bt}}},\O_Z),
$$
for $i: Z \to Z^{\Delta^n}$. However, $Z$ is quasi-\'etale (the analogous notion to quasi-smooth), so the vertex maps $a:Z^{\Delta^n} \to Z$ are trivially \'etale, and thus $\Omega_{Z^{\Delta^n}} \cong a^*\Omega_Z $, so $i^*\Omega_{Z^{\Delta^{\bt}}}= \Omega_Z$. Therefore $\H^0E=\H^0E_0=\Hom_{\O_X}(\Omega_X, \O_X)$, and $\H^nE=0$ for $n<0$.
\end{proof}

\subsection{Quasi-compact, quasi-separated stacks}\label{nicestacks}

Let $\fX$ be a quasi-compact, quasi-separated stack, with presentation $P:X_0 \to \fX$, for $X_0$ affine, giving a  simplicial algebraic space $\cosk_0^{\fX}(X_0)$ (as considered in \cite{aoki} \S 3). We may then take an \'etale hypercovering $X_{\bt} \to \cosk_0^{\fX}(X_0)$, for $X_{\bt}$ a simplicial affine scheme, and denote the composition by $P_{\bt}:X_{\bt} \to \fX$. 

\begin{lemma}
Every smooth simplicial hypercovering is trivially smooth.
\end{lemma}
\begin{proof}
For a map $U_{\bt} \to V_{\bt}$ to be a smooth hypercovering says that the matching maps $U_n \to V_n\by_{M_nV}M_nU$ are all smooth surjections.
\end{proof}

\begin{lemma}
The simplicial affine scheme $X_{\bt}$ is quasi-smooth.
\end{lemma}
\begin{proof}
Write $Z_{\bt}:= \cosk_0^{\fX}(X_0)$.
Since $Z_{\bt}=B\fG$, for $\fG$ the groupoid space $\xymatrix@1{X_0\by_{\fX}X_0 \ar@<.5ex>[r] \ar@<-.5ex>[r] & X_0}$, all higher partial matching maps of $Z_{\bt}$ are isomorphisms. In other words, for any trivial cofibration $i:K \to L$ in $\bS$ with $i_0:K_0 \to L_0$ an isomorphism, the map
$$
i^*: M_L(Z) \to M_K(Z)
$$
is an isomorphism.

By \cite{aoki} Theorem 2.1.5, $\fG$  has SQCS structure so  the maps  $X_1 \to \Hom_{\bS}(\L^1_k,X_{\bt})$ are smooth surjections for both $k$. Thus $Z_{\bt}$ is quasi-smooth, by Lemma \ref{gsmoothslack}. Since $X_{\bt} \to Z_{\bt}$ is trivially smooth, the result follows.
\end{proof}

\begin{remark}\label{inftyrk}
Similarly, every strongly quasi-compact $n$-geometric Artin stack  $\fX$    gives rise to a quasi-smooth simplicial affine scheme $X_{\bt}$, by \cite{stacks2} Theorem \ref{stacks-relstrict}. The statement of Proposition  \ref{olssoncoho} will then carry over to this generality, taking $\bL_{\fX}$ to be the cotangent complex of \cite{stacks2} \S \ref{stacks-cotsn}.
\end{remark}

\subsubsection{Cohomology and the cotangent complex}\label{cfolsson}

Given any morphism $f:\fY \to \fX$ of quasi-compact, quasi-separated stacks, lifting to a morphism $f:Y_{\bt} \to X_{\bt}$ of simplicial affine resolutions,  in this section we will describe the $\Ext$-groups 
$$
\Ext^*_{c\Alg/c_+\Mod}(f^{\sharp})
$$
of Proposition \ref{stackqsgood} in terms of the cotangent complex of \cite{olssartin} \S 8. 
 $\Ext$-groups of the cotangent complex are defined in \cite{olssonstack} \S 2.11.

Let $\fX, X_{\bt}$ be as above, and let $\sJ$ be a quasi-coherent sheaf on $\fX$. Since the cotangent complex $\bL_{\fX}$ is in degrees $\ge -1$, we have $\Ext^i(\bL_{\fX}, \sJ)=0$ for all $i <-1$. Since $r:X_{\bt} \to \cosk_0^{\fX}(X_0)$ is a hypercovering, the maps
$$
\H^*(\cosk_0^{\fX}(X_0), \sF) \to \H^*(X_{\bt}, r^*\sF)
$$
on cohomology  are isomorphisms for all quasi-coherent sheaves $\sF$.

By \cite{aoki} Proposition 3.4.2, 
\begin{enumerate}
\item
$\Ext^i(\bL_{\fX}, \sJ)\cong \Ext^i(\bL_{X_{\bt}},P_{\bt}^*\sJ)$  for $i>0$. 

\item $\Ext^i(\bL_{\fX}, \sJ)= \H^i( \Hom(\Omega_{X_0/\fX}, P^*\sJ)[1] \to \Hom(\Omega_{X_{\bt}}, P_{\bt}^*\sJ))$, for $i \le 0$.
\end{enumerate}

\begin{proposition}\label{olssoncoho}
$$
\Ext^i_{c\Alg/c_+\Mod}(f^{\sharp})\cong \Ext^i(\bL_{\fX}, f_*\O_{\fY}) 
$$
for all $i \in \Z$.
\end{proposition}
\begin{proof}
For $i>0$, this is just the observation that $\Ext^i(\bL_{X_{\bt}},P_{\bt}^*\sJ)=\Ext^i_{c\Alg/c_+\Mod}(f^{\sharp})$ when $\sJ=f_*\O_{\fY}$.

Accordingly, we need to describe the non-positive $\Ext$ groups
$$
\H_{-i}\Hom_{X}(c^*\Omega_{X^{\Delta^{\bt}}},\sJ)
$$
in terms of $\fX, X_0$, where $c: X \to X^K$ is the constant map.

Let $U_{\bt}$ denote the simplicial complex  $\Hom_{X}(c^*\Omega_{X^{\Delta^{\bt}}},\sJ)$, and 
write $Z_{\bt}:= \cosk_0^{\fX}(X_0)$, with $V_{\bt}: \Hom_{Z}(c^*\Omega_{Z^{\Delta^{\bt}}},\sJ)$. Since $X \to Z$ is trivially smooth, observe that the canonical map $U \to V$ is a trivial fibration, so $\H_*(U) \cong \H_*(V)$.

In general, if $K$ is contractible, then 
$$
M_KZ= \overbrace{X_0\by_{\fX}X_0\by_{\fX}\ldots \by_{\fX}X_0}^{K_0},
$$
so 
$$
\Omega(M_KZ/\fX)= \bigoplus_{v \in K_0} v^*\Omega(X_0/\fX).
$$

We therefore conclude that for a trivial cofibration $K \into L$, 
$$
\Omega(M_LZ/M_KZ)= \bigoplus_{v \in L_0 - K_0} v^*\Omega(X_0/\fX),
$$
so for $0:\bt \to I$, 
$$
\Omega(Z^I/Z)^n= \Omega(M_{I\by \Delta^n}Z/M_{\Delta^n}Z) =\bigoplus_{v \in \Delta^n_0} (v \by 1)^*\Omega(X_0/\fX),
$$
so 
$$
i^*\Omega(Z^I/Z)^n= \bigoplus_{v \in \Delta^n_0} v^*\Omega(X_0/\fX), 
$$
and 
$$
\Hom_{\O(X)}(i^*\Omega(X^I/X), \sF)= \Hom_{\O(X_0)}(\Omega(X_0/\fX), \sF^0),
$$
giving
$$
N_1V= \Hom_{\O(X)}(i^*\Omega(X^I/X),f_*\O(Y))= \Hom_{\O(X_0)}(\Omega(X_0/\fX), f_*\O(Y_0)).
$$

Moreover, for $n \ge 2$, $\L^n_0 \to \Delta^n_0$ is an isomorphism, so $N_nV=0$.

Thus 
$$
NV= (\Hom(\Omega_{X_0/\fX}, P^*\sJ)[-1] \to \Hom(\Omega_{X_{\bt}}, P_{\bt}^*\sJ)),
$$
as required.
\end{proof}

\subsubsection{Comparing deformation groupoids}\label{cfaoki}

\begin{definition}
Given a small $2$-category $\C$, define a simplicial category $B^1\C$ by setting $\Ob(B^1\C)= \Ob(\C)$, and $\HHom_{B^1\C}(x,y)= B\hom_{\C}(x,y)$, where $\hom_{\C}(x,y)$ is the $1$-category of homomorphisms from $x$ to $y$, and $B$ is the nerve functor. 
\end{definition}

\begin{lemma}\label{pi2gpd}
Given $x,y \in \Ob \C$   with $\hom_{\C}(x,y)$ a groupoid,   $\pi_0\HHom_{B^1\C}(x,y)$ is the set of isomorphism classes in $\hom_{\C}(x,y)$, with $\pi_1(\HHom_{B^1\C}(x,y),f)$ the set of $2$-automorphisms of $f$, and $\pi_i(\HHom_{B^1\C}(x,y),f)=0$ for $i>1$.
\end{lemma}

\begin{definition}\label{2catgpd}
Define a $2$-category structure on the category $\Alg\gpd\Sp$ of algebraic groupoid spaces  (as in \cite{aoki}) by defining a $2$-morphism $\eta $ between morphisms $f,f':G \to H$ by analogy with natural transformations. Explicitly, let $\Ob G$ be the space of objects of $G$, with $\Mor G \to (\Ob G \by \Ob G)$ the space of isomorphisms, and similarly for $H$. We must have $\eta:\Ob G \to \Mor H$, with $s\circ \eta=f, t \circ \eta= f'$, and the following diagram commuting
$$
\begin{CD}
\Mor (G) @>{(\eta \circ t, \Mor f)}>> \Mor (H)\by_{s, \Ob H, t}\Mor(H) \\
@V{(f', \eta \circ s) }VV @VVmV\\
\Mor (H)\by_{s, \Ob H, t}\Mor(H) @>m>> \Mor (H).
\end{CD}
$$
  \end{definition}

\begin{definition}
Given a $2$-groupoid $\cG$, define $\Pi_0\cG$ to be the groupoid with objects $\Ob \cG$, and morphisms $\Hom_{\Pi_0\cG}(X,Y)= \pi_0\hom_{\cG}(X,Y)$. Similarly, for  a simplicial groupoid $G_{\bt}$, define $\Pi_0G_{\bt}$ to be the groupoid with objects $\Ob G$, and morphisms $\Hom_{\Pi_0\G_{\bt}}(X,Y)= \pi_0\HHom_{G}(X,Y)$.
\end{definition}

\begin{lemma}
Given $G \in \Alg\gpd\Sp$ associated to an algebraic stack over $k$,  the nerve functor $B: \Alg\gpd\Sp\to s\Alg\Sp$ to the category of simplicial algebraic spaces gives an isomorphism
$$
B^1\Def^2_{\Alg\gpd\Sp}(G) \cong \underline{\Def}_{s\Alg\Sp}(BG),
$$
between the  $2$-groupoid of deformations in $\Alg\gpd\Sp$, and the simplicial groupoid of deformations in $s\Alg\Sp$.
\end{lemma}
\begin{proof}
By \cite{aoki} Corollary 3.1.5, we know that $\Pi_0\Def^2_{\Alg\gpd\Sp}(G) \cong \Pi_0\underline{\Def}_{s\Alg\Sp}(BG)$, so  we just need to show that, for algebraic groupoid spaces $H,G$, 
$$
\HHom_{s\Alg\Sp}(BH,BG)=B\hom_{\Alg\gpd\Sp}(H,G).
$$
Now, $\HHom_{s\Alg\Sp}(X,BG)_n= \Hom_{s\Alg\Sp}(X\by \Delta^n, BG)= \Hom_{\Alg\gpd\Sp}(\pi_fX \by \pi_f\Delta^n, G)$, where we define the fundamental groupoid $\pi_f:s\Alg\Sp\to \Alg\gpd\Sp$ to be left adjoint to $B$, noting that $\pi_fBH=H$. However, $\pi_f\Delta^n$ is the groupoid with $n+1$ objects, and unique isomorphisms between them. Thus
$$
\Hom_{\Alg\gpd\Sp}(\pi_fX \by \pi_f\Delta^n, G)= B_n \hom_{\Alg\gpd\Sp}(\pi_fX,G),
$$
as required.
\end{proof}

\begin{lemma}
The functor $C$ defined in \cite{aoki} \S3.2 gives an equivalence between  $ \Def^2_{\Alg\gpd\Sp}(G)$ and $\Def^2(CG)$, the $2$-groupoid of deformations of the algebraic stack $CG$.
\end{lemma}
\begin{proof}
First observe that $C$ maps the $2$-isomorphisms of Definition \ref{2catgpd} to $2$-isomorphisms of stacks, so $C$ is well-defined. 

By \cite{aoki} Proposition 3.2.5, we know that $C$ induces a bijection on isomorphism classes of objects. 
We need to show that for $\cG \in \Def^2_{\Alg\gpd\Sp}(G)(A)$, 
$$
C:\hom_{\Def^2_{\Alg\gpd\Sp}(G)(A)}(\cG,\cG)\to \hom_{\Def^2(CG)(A)}(C\cG,C\cG)
$$
 is an equivalence of groupoids. By [ibid.] Proposition 3.3.2, it is essentially surjective. 

Given $f \in \Ob\hom_{\Def^2_{\Alg\gpd\Sp}(G)(A)}(\cG,\cG)$, we thus need to show that 
$$
\theta:\Aut^2_{\Def^2_{\Alg\gpd\Sp}(\cG)(A)}(f) \to \Aut^2_{ \Def^2(C\cG)(A)}(Cf)
$$
 is an isomorphism of $2$-automorphism groups. Multiplication by $f^{-1}$ allows us to assume that $f=\id_{\cG}$.

By [ibid.] Proposition 3.3.2, we 
have an exact sequence
$$
0 \to \Aut^2_{ \Def^2(C\cG)}(\id_{C\cG}) \to \Aut(X_0/C\cG)_P \xra{A} \Aut_{\Def^2_{\Alg\gpd\Sp}(G)}(\cG).
$$
Since $\Aut(X_0/C\cG)_P$ is smooth, the homogeneous functor $\Aut^2_{ \Def^2(C\cG)}(\id_{C\cG})$ has tangent space $\ker(\tan A)$ and obstruction space $\coker(\tan A)$. By [ibid.] Proposition 3.4.2, these are  $\Ext^{-1}(\bL_{\fX}, \O_{fX})$ and  $\Ext^{0}(\bL_{\fX}, \O_{fX})$, respectively. Thus \S \ref{cfolsson}, Lemma \ref{pi2gpd} and Theorem \ref{robs} imply that $\theta$  gives isomorphisms on tangent and obstruction spaces, so must be an isomorphism of homogeneous functors by the standard smoothness criterion.
\end{proof}

Taking $G= \xymatrix@1{X_0\by_{\fX}X_0 \ar@<.5ex>[r] \ar@<-.5ex>[r] & X_0}$, we have therefore shown that the deformation $2$-groupoid of $\fX$ is equivalent to the simplicial deformation groupoid of $\cosk^{\fX}_0(X_0)$. We still need to compare this with the simplicial affine scheme $X_{\bt}$ defined at the beginning of the section.

\begin{proposition}
The simplicial deformation groupoids of $\cosk^{\fX}_0(X_0)$ and $X_{\bt}$ are equivalent.
\end{proposition}
\begin{proof}
Let $Z_{\bt}:=\cosk^{\fX}_0(X_0)$. As in \S \ref{weakinvart}, we will consider three simplicial deformation problems $F_X, F_Z, F_r$: deformations of $X_{\bt}$,  deformations of $Z_{\bt}$, and deformations of the diagram $r:X_{\bt} \to Z_{\bt}$. Note that these all define quasi-smooth functors $F:\C_{\L} \to s\gpd \xra{\bar{W}} \bS$, so we just need to compare tangent and obstruction spaces. 

The calculations of Proposition \ref{olssoncoho} show that the loop spaces  $\Omega F_X, \Omega F_Z$ have tangent spaces
$$
\underline{\Hom}_{\O_X}(\Omega_X, \O_X), \quad \underline{\Hom}_{\O_Z}(\Omega_Z, \O_Z).
$$
Similarly, there is a fibration $ F_r \to  F_X \by  F_Z$, whose fibre has tangent space 
$$
\underline{\Hom}_{\O_Z}(\Omega_Z, r_*\O_X).
$$
Since the maps
$$
\underline{\Hom}_{\O_X}(\Omega_X, \O_X) \xra{r_*} \underline{\Hom}_{\O_Z}(\Omega_Z, r_*\O_X) \xla{r^*}\underline{\Hom}_{\O_Z}(\Omega_Z, \O_Z)
$$
are isomorphisms, we deduce that the maps $ F_Z \la F_r \to F_X $ induce isomorphisms on tangent spaces of positive homotopy groups.

It only remains to show that the deformation functors $\pi_0F_X, \pi_0F_Z, \pi_0F_r $ have isomorphic tangent and obstruction spaces. By adapting \cite{paper2} \S \ref{paper2-etshf}, we may deduce that these are (respectively)
$$
\Ext^i_{\O_X}(\bL_{X}, \O_X), \quad \Ext^i_{\O_Z}(\bL_{Z}, \O_Z), 
$$
and the groups $T^i$ fitting into the long exact sequence
$$
\ldots \to T^1 \to \Ext^1_{\O_X}(\bL_{X}, \O_X)\by \Ext^1_{\O_Z}(\bL_{Z}, \O_Z)\to \Ext^1_{\O_Z}(\bL_{Z}, r_*\O_X) \to T^2 \to \ldots .
$$
Now, since the maps
$$
\Ext^i_{\O_X}(\bL_{X}, \O_X)\xra{r_*}\Ext^i_{\O_Z}(\bL_{Z}, r_*\O_X) \xla{r^*} \Ext^1_{\O_Z}(\bL_{Z}, \O_Z)
$$
are isomorphisms for $i\ge 1$, with $r_*$ surjective for $i=0$, we see that the functors $F_r,F_X,F_Z$ are all equivalent.
\end{proof}

\subsection{Arbitrary algebraic stacks}

We now wish to describe derived deformations of  a simplicial scheme $X_{\bt}$ over $k$, with each $X_n$ a disjoint union of affine schemes. 

\begin{definition}
For any scheme $Y$, let $\pi(Y)$ be the set of connected components of $Y$, and $\pi: Y \to \pi(Y)$ the map of associated topological spaces.
\end{definition}

Now, deformations of $X_n$ are equivalent to deformations of the algebra $\pi_*\O_{X_n}$ over $\pi(X_n)$.

\begin{definition}
Recall that the ordinal number categories $\Delta, \Delta^*$ can be regarded as subcategories of $\bS$, by identifying $[n]$ with $\Delta^n$.

Given a category $\C$ and $K \in \bS$, define $c\C^K$ (resp. $c_+\C^K$) to be the category of functors from $\Delta \da K$ (resp. $\Delta_* \da K$) to $\C$. Thus an object $C \in c\C^K$ consists of objects $C_a$ for all $n \in \N_0$, $a \in K_n$, together with compatible maps $\pd^i: M_{\pd_ia} \to M_a$, $\sigma^i:  M_{\sigma_ia} \to M_a$, and similarly for $c_+\C^K$.
\end{definition}

Now, observe that $\pi_*\O_X$ defines an object of $c\Alg(k)^{\pi(X)}$, with 
$$
(\pi_*\O_X)_a=\Gamma (\pi^{-1}(a), \O_{X_n}),
$$
for $a \in \pi(X_n)$. Since any deformation of $X_{\bt}$ will not change $\pi(X)$, deformations of $X_{\bt}$ are equivalent to deformations of $\pi_*\O_X$.

The categories $c\FAlg^{\pi(X)}, c_+\FAlg^{\pi(X)}, c\FMod^{\pi(X)}, c_+\FMod^{\pi(X)}$ can all be given simplicial structures as in Definition \ref{csimpstr}, setting $(C^K)_a= (C_a)^{K_n}$ for $a \in \pi(X_n)$.

\begin{remark}
Observe that for any category $\C$ and any map $f:K \to L$ in $\bS$, there are maps $f^{-1}: c\C^L \to  c\C^K$, $f^{-1}: c_+\C^L \to  c_+\C^K$ given by $(f^{-1}C)_a=C_{f(a)}$. If $\C$ contains products, then $f^{-1}$ has a right adjoint $f_*$, given by $(f_*C)_b= \prod_{a \in f^{-1}(b)}C_a$. For $f:K\by L \to L$, $C^K= f_*f^{-1}C$.
\end{remark}

If $f:\pi(X) \to \bt$ denotes the constant map, then we write $\Gamma:= f_*$, with the constant functor $f^{-1}$ denoted by $\Gamma^*$.

We then have a diagram
 of adjunctions of functors $\C_{\L} \to s\Cat$:
$$
\xymatrix@=8ex{
c\FAlg(A)^{\pi(X)} \ar@<1ex>[r]^-{U_{\pd}U_{\alg}}_-{\top} \ar@<-1ex>[d]_{\Gamma}^{\vdash} 
&\ar@<1ex>[l]^-{ \Symm  F_{\pd}}   c_+\FMod(A)^{\pi(X)}    \ar@<-1ex>[d]_{\Gamma}^{\vdash} 
\\
\ar@<-1ex>[u]_{\Gamma^*}c\FAlg(A)	 \ar@<1ex>[r]^-{U_{\pd}U_{\alg}}_-{\top} 
&\ar@<1ex>[l]^-{\Symm F_{\pd} } \ar@<-1ex>[u]_{\Gamma^*}  c_+\FMod(A), 
}
$$
 where $F_{\pd}:c_+\C^{\pi(X)}\to c\C^{\pi(X)}$ is left adjoint to the forgetful functor $U_{\pd}$, given by 
$$
(F_{\pd} C^*)_a = C_a\sqcup C_{\pd_0a} \sqcup \ldots \sqcup C_{(\pd_0)^n a},
$$
for $a \in \pi(X_n)$, with operations dual to those in Lemma \ref{Gdef}. 

We must check that $\Gamma^* \dashv \Gamma$ is monadic.
 For this, we verify Beck's Theorem (e.g. \cite{Mac} Ch. VI.7 Ex. 6), observing that $\Gamma$ commutes with coequalisers --- this is effectively the observation that taking arbitrary products is an exact functor. 

Writing $U:=U_{\pd}U_{\alg}$ and $F:= \Symm F_{\pd}$,  we also have the following commutativity conditions:
$$
\Gamma U=U \Gamma \quad \Gamma^*F=F \Gamma^*, \quad  \Gamma^* U= U \Gamma^*,
$$
and a natural transformation
$$
F\Gamma\to \Gamma F.
$$

These adjunctions combine to give a monadic adjunction
$$
\xymatrix@1{ c\FAlg(A)^{\pi(X)} \ar@<1ex>[r]^-{\Gamma U_{\pd} U_{\alg}}_-{\top}  
&\ar@<1ex>[l]^-{\Gamma^* \Symm  F_{\pd}} c_+\FMod(A)}.
$$

\begin{definition}
Given a simplicial scheme $X_{\bt}$, with each $X_n$ a disjoint union of affine schemes, define $c\Mod(X)$ to be the category of $\pi_*(\O_X)$-modules over $c\Mod^{\pi(X)}$.
\end{definition}

\begin{lemma}\label{omegaproj2}
If $X \to Y$ is  a trivially smooth map of simplicial schemes, with each $X_n,Y_n$ a disjoint union of affine schemes, then $\pi_*\Omega(X/Y)$ is projective in $c\Mod(X)$.
\end{lemma}
\begin{proof}
This is similar to Lemma \ref{omegaproj}. We may define matching objects of  $L \in c\Mod(X)$ by letting $M^nL$ on $\pi(X_{n+1})$ be the equaliser 
$$
\xymatrix@1{M^nL \ar[r] & \prod_{i=0}^{n} \sigma_{i*}L^{n} \ar@<.5ex>[r]^-{a} \ar@<-.5ex>[r]_-b & \prod_{0 \le i <j \le n} \sigma_{i*}\sigma_{j*}L^{n-1}}, 
$$
where
$\pr_{ij} \circ a= \sigma^i\circ \pr_j$, $\pr_{ij} \circ b= \sigma^{j-1} \circ \pr_i$.
Note that  $\Gamma(M^nL)=M^n(\Gamma L)$. Since $\Gamma$ reflects isomorphisms, this means that for all surjections $L \onto N$, the relative matching map $L^n \to M^{n-1}L\by_{M^{n-1}N}N^n$ is surjective.

In order to construct latching maps, note that any cocontinuous functor $S:(\Delta \da \pi(X)) \to \Alg$ extends to a cocontinuous functor $S:(\bS \da \pi(X)) \to \Alg$.
Given   $M^{\bt}  \in c\Mod(X)$ and $a:K\to \pi(X) $ in $\bS$, define $M(a) \in \Mod(\pi_*(\O_X)(a))$ by 
$$
\pi_*(\O_X)(a) \oplus M(a)\eps=(\pi_*(\O_X)\oplus M\eps)(a).
$$ 
Note that if we set $X(a):= \Hom_{\bS\da \pi(X)}(K, X)$, then   $\Omega(X/Y)(a)=\Omega(X(a)/Y(\pi(f)_*a))$.

The latching object of $\Omega(a)$, for $a \in \pi(X)_n$, is $\Omega(\pd a)$, for $\pd:\pd\Delta^n \to \Delta^n$. It therefore suffices to show that $(X^{\pd})^*\Omega(\pd a)\to \Omega(a)$ is  projective in $\Mod(X(a))$ for all such $a$. By adapting the proof of Lemma \ref{omegaproj}, it suffices to show that
$$
X(a) \to Y(\pi(f)_*a)\by_{Y(\pi(f)_*\pd a)}X(\pd a)
$$
is smooth.

Set $Y':=X\by_{\pi(X)}\pi(Y)$, and observe that Lemma \ref{smoothcdns} implies that $X \to Y'$ is trivially smooth. Thus the matching map $X_n \to M_nX\by_{M_nY'}Y_n'$ is smooth. The required result is then obtained by taking the fibre over $ a \in \pi(X)_n$.
\end{proof}

\begin{lemma}\label{lalg2}
If we set $\bot'= \Gamma^* FU\Gamma$, and  
$$
\bL^{\bot'}_{\bt}(X):= \Omega((\bot')^{n+1}\pi_*(\O_X))\ten_{(\bot')^{n+1}\pi_*(\O_X)} \pi_*(\O_X),
$$
 then 
for all $m$, the simplicial complex $\bL^{\bot'}_{\bt}(X)^m$ is a model for the cotangent complex of $X_m$.
\end{lemma}
\begin{proof}
This is essentially the same as Lemma \ref{lalg}, making use of the observation that $\Gamma$ is exact and reflects isomorphisms, so it suffices to prove that $U\Gamma(\bot')^{\bt +1}\pi_*(\O_X)\to U\Gamma\pi_*(\O_X)$ is a resolution.
\end{proof}

\begin{definition}
Define $\cD$ to be the simplicial category of pairs $(K, R)$, for $K \in \bS$, $R \in (c\FAlg^K)^{\op}$, with a  morphism $f \in \HHom_{\cD}( (K,R), (L, S))_n$ consisting of $f:K \to L$ in $\bS$, together with $f^{\sharp} \in \HHom_{c\FAlg^K}( f^{-1}S, R)_n$. 

Define $ \cB:= \bS \by (c_+\Mod)^{\op} $, with simplicial structure coming from $(c_+\Mod)^{\op}$.
\end{definition}

Now, observe that we have a forgetful functor $V:\cD \to \cB$, given by $(K, R) \mapsto (K, \Gamma U_{\pd} U_{\alg} R)$, with right adjoint $G:\cB\to \cD$ given by $(K, M) \mapsto (K,\Gamma^* \Symm  F_{\pd}M)$. We have already seen that this adjunction is comonadic (by fixing $K$). 

\begin{proposition}\label{stackqsgood2}
 If $X,Y$  are  simplicial schemes over $k$, with each $X_n, Y_n$ a disjoint union of affine schemes, and $X$ quasi-smooth, then every   morphism $\rho:Y \to X$   is quasi-smooth over $\cB$, in the sense of Definition \ref{qsmorphism}.
\end{proposition}
\begin{proof}
The proof of Proposition \ref{stackqsgood} carries over, using Lemmas \ref{omegaproj2} and \ref{lalg2} instead of  Lemmas \ref{omegaproj} and \ref{lalg}.
\end{proof}

\begin{corollary}
For any diagram in $\cD(k)$ with quasi-smooth objects, the associated pre-DDC given by Definition \ref{sdcdiagram} and Proposition \ref{enrich} applied to the adjunction $G \vdash V$ 
is a DDC by Lemma \ref{q2diagram}, and governs deformations in the simplicial category $\cD$ (by Proposition \ref{governdiagram}).
\end{corollary}

\bibliographystyle{alphanum}
\addcontentsline{toc}{section}{Bibliography}
\bibliography{references}
\end{document}

%% file: newthms.tex
\newtheorem{theorem}{Theorem}[section]
\newtheorem{proposition}[theorem]{Proposition}
\newtheorem{corollary}[theorem]{Corollary}

\newtheorem{lemma}[theorem]{Lemma}
\newtheorem*{theorem*}{Theorem}
\newtheorem*{proposition*}{Proposition}
\newtheorem*{corollary*}{Corollary}
\newtheorem*{lemma*}{Lemma}
\newtheorem*{conjecture*}{Conjecture}

\theoremstyle{definition}
\newtheorem{definition}[theorem]{Definition}

\newtheorem*{definition*}{Definition}

\theoremstyle{remark}
\newtheorem{example}[theorem]{Example}
\newtheorem{examples}[theorem]{Examples}
\newtheorem{remark}[theorem]{Remark}
\newtheorem{remarks}[theorem]{Remarks}

\newtheorem*{example*}{Example}
\newtheorem*{examples*}{Examples}
\newtheorem*{remark*}{Remark}
\newtheorem*{remarks*}{Remarks}
\newtheorem*{exercise*}{Exercise}

%% file: newcommands.tex
\newcommand\da{\!\downarrow\!}

\newcommand\la{\leftarrow}

\newcommand\id{\mathrm{id}}

\newcommand\ten{\otimes}
\newcommand\vareps{\varepsilon}
\newcommand\eps{\epsilon}

\newcommand\CC{\mathrm{C}}

\renewcommand\H{\mathrm{H}}
\newcommand\z{\mathrm{Z}}

\newcommand\N{\mathbb{N}}
\newcommand\Z{\mathbb{Z}}

\newcommand\bD{\mathbb{D}}

\newcommand\bI{\mathbb{I}}
\newcommand\bJ{\mathbb{J}}

\newcommand\bL{\mathbb{L}}

\newcommand\bS{\mathbb{S}}

\newcommand\C{\mathcal{C}}

\newcommand\cA{\mathcal{A}}
\newcommand\cB{\mathcal{B}}

\newcommand\cD{\mathcal{D}}
\newcommand\cE{\mathcal{E}}
\newcommand\cF{\mathcal{F}}
\newcommand\cG{\mathcal{G}}

\newcommand\cL{\mathcal{L}}

\newcommand\cO{\mathcal{O}}

\newcommand\cS{\mathcal{S}}

\newcommand\G{\mathscr{G}}

\renewcommand\O{\mathscr{O}}

\newcommand\sF{\mathscr{F}}
\newcommand\sG{\mathscr{G}}

\newcommand\sJ{\mathscr{J}}

\newcommand\sM{\mathscr{M}}
\newcommand\sN{\mathscr{N}}

\newcommand\Def{\mathfrak{Def}}

\newcommand\fG{\mathfrak{G}}

\newcommand\fX{\mathfrak{X}}
\newcommand\fY{\mathfrak{Y}}

\renewcommand\L{\Lambda}

\newcommand\m{\mathfrak{m}}

\renewcommand\hom{\mathscr{H}\!\mathit{om}}

\newcommand\Ho{\mathrm{Ho}}
\newcommand\Ring{\mathrm{Ring}}
\newcommand\Alg{\mathrm{Alg}}
\newcommand\Mod{\mathrm{Mod}}
\newcommand\FAlg{\mathrm{FAlg}}
\newcommand\FMod{\mathrm{FMod}}
\newcommand\Hom{\mathrm{Hom}}

\newcommand\HHom{\underline{\mathrm{Hom}}}

\newcommand\Ext{\mathrm{Ext}}
\newcommand\EExt{\mathbb{E}\mathrm{xt}}

\newcommand\Der{\mathrm{Der}}

\newcommand\Aut{\mathrm{Aut}}

\newcommand\Iso{\mathrm{Iso}}

\newcommand\coker{\mathrm{coker\,}}

\newcommand\Ob{\mathrm{Ob}\,}

\newcommand\Gp{\mathrm{Gp}}

\newcommand\Spec{\mathrm{Spec}\,}

\newcommand\Spf{\mathrm{Spf}\,}

\newcommand\Set{\mathrm{Set}}

\newcommand\Cat{\mathrm{Cat}}

\newcommand\Dat{\mathrm{Dat}}

\newcommand\Mor{\mathrm{Mor}\,}

\newcommand\Sp{\mathrm{Sp}}

\newcommand\ad{\mathrm{ad}}

\newcommand\into{\hookrightarrow}
\newcommand\onto{\twoheadrightarrow}
\newcommand\abuts{\implies}
\newcommand\xra{\xrightarrow}
\newcommand\xla{\xleftarrow}
\newcommand\pr{\mathrm{pr}}

\newcommand\alg{\mathrm{alg}}

\newcommand\bt{\bullet}
\newcommand\by{\times}

\newcommand\mc{\mathrm{MC}}
\newcommand\mmc{\underline{\mathrm{MC}}}
\newcommand\MC{\mathfrak{MC}}

\newcommand\ddef{\mathrm{Def}}

\newcommand\Vect{\mathrm{Vect}}

\newcommand\Symm{\mathrm{Symm}}

\newcommand\Tot{\mathrm{Tot}\,}
\newcommand\diag{\mathrm{diag}\,}
\newcommand\toph{\top_{\mathrm{h}}}
\newcommand\topv{\top_{\mathrm{v}}}
\newcommand\both{\bot_{\mathrm{h}}}
\newcommand\botv{\bot_{\mathrm{v}}}

\newcommand\pd{\partial}

\newcommand\gpd{\mathrm{Gpd}}

\renewcommand\alg{\mathrm{alg}}

\newcommand\cosk{\mathrm{cosk}}

\newcommand\op{\mathrm{opp}}

\newcommand\oR{\mathbf{R}}

\newcommand\uleft\underleftarrow
\newcommand\uline\underline